\documentclass[12pt] {article}
\usepackage{amsmath,amsthm}
\usepackage{amssymb,latexsym}
\usepackage{amscd}
\title{Theory of valuations on manifolds, IV. New properties of the multiplicative structure.}
\date{}
\author{ Semyon Alesker \footnote{Partially supported by ISF grant 1369/04.}
\\  { \normalsize Department of Mathematics, Tel Aviv University, Ramat Aviv}
 \\  { \normalsize 69978 Tel Aviv,
Israel }
\\ {\normalsize e-mail: semyon@post.tau.ac.il}}

\def\eps{\varepsilon}
\def\alp{\alpha}
\def\ome{\omega}
\def\Ome{\Omega}
\def\lam{\lambda}
\def\Lam{\Lambda}

\def\to{\rightarrow}
\def\qed { Q.E.D. }
\def\pt{\partial}

\def\si{{}\!^{S}\int}
\def\sio{{}\!^ {S_1}\int}
\def\sit{{}\!^ {S_2}\int}

\def\RR{\mathbb{R}}
\def\CC{\mathbb{C}}

\def\NN{\mathbb{N}}
\def\ZZ{\mathbb{Z}}

\def\DD{\mathbb{D}}
\def\PP{\mathbb{P}}

\def\One{{1\hskip-2.5pt{\rm l}}}
\swapnumbers
\newtheorem{theorem}{Theorem}[subsection]
\newtheorem{corollary}[theorem]{Corollary}
\newtheorem{lemma}[theorem]{Lemma}
\newtheorem{proposition}[theorem]{Proposition}

\theoremstyle{definition}

\newtheorem{definition}[theorem]{Definition}
\newtheorem{remark}[theorem]{Remark}

\theoremstyle{proposition-definition}
\newtheorem{proposition-definition}[theorem]{Proposition-Definition}

\numberwithin{equation}{subsection}

\def\cf{{\cal F}}

  \def\cf{{\cal F}}
  
 \def\ck{{\cal K}} \def\cl{{\cal L}}
  \def\co{{\cal O}}
\def\cp{{\cal P}} 

 \def\ct{{\cal T}} 
\def\cv{{\cal V}} \def\cw{{\cal W}} 
 
\def\mov{|\omega_V|}
\def\inj{\hookrightarrow }

\def\vi{V^\infty(X)}
\def\vmi{V^{-\infty}(X)}
\def\vic{V^\infty_c(X)}
\def\vmic{V^{-\infty}_c(X)}
\newcommand \supp{\operatorname{supp} \,}
\begin{document}
\maketitle \setcounter{section}{-1}
\begin{abstract}
This is the fourth part in the series of articles \cite{part1},
\cite{part2}, \cite{part3} (see also \cite{alesker-poly}) where
the theory of valuations on manifolds is developed. In this part
it is shown that the filtration on valuations is compatible with
the product. Then it is proved that the Euler-Verdier involution
on smooth valuations is an automorphism of the algebra of
valuations. Then an integration functional on valuations with
compact support is introduced, and a property of selfduality of
valuations is proved. Next a space of generalized valuations is
defined, and some basic properties of it are proved. Finally a
canonical imbedding of the space of constructible functions on a
real analytic manifold into the space of generalized valuations is
constructed, and various structures on valuations are compared
with known structures on constructible functions.
\end{abstract}
\tableofcontents
\section{Introduction.}
\setcounter{subsection}{1} In convexity there are many
geometrically interesting and well known examples of valuations on
convex sets: Lebesgue measure, the Euler characteristic, the
surface area, mixed volumes, the affine surface area. For a
description of older classical developments on this subject we
refer to the surveys \cite{mcmullen-schneider},
\cite{mcmullen-survey}. For the general background on convexity we
refer to the book \cite{schneider-book}.

Approximately during the last decade there was a significant
progress in this classical subject which has led to new
classification results of various classes of valuations, to
discovery of new structures on them. This progress has shed a new
light on the notion of valuation which allowed to generalize it to
more general setting of valuations on manifolds and on not
necessarily convex sets (which do not make sense on a general
manifold). On the other hand author's feeling is that the notion
of valuation equips smooth manifolds with  a new general
structure. The development of the theory of valuations on
manifolds was started in three previous parts of the series of
articles: \cite{part1}, \cite{part2} by the author and
\cite{part3} by J. Fu and the author. This article in the forth
part in this series.

 In \cite{part2} the notion of smooth
valuation on a smooth manifold was introduced. Roughly speaking a
smooth valuation can be thought as a finitely additive
$\CC$-valued measure on a class of nice subsets; this measure is
requested to satisfy some additional assumptions of continuity (or
rather smoothness) in some sense. The basic examples of smooth
valuations on a general manifold $X$ are smooth measures on $X$
and the Euler characteristic. Moreover the well known intrinsic
volumes of convex sets can be generalized to provide examples of
smooth valuations on an arbitrary {\itshape Riemannian} manifold;
these valuations are known as Lipschitz-Killing curvatures.

Let $X$ be a smooth manifold of dimension $n$. The space of smooth
valuations on $X$ is denoted by $V^\infty(X)$. It has a canonical
linear topology with respect to which it becomes a Fr\'echet
space.

The space $V^\infty(X)$ carries a canonical multiplicative
structure. This structure seems to be of particular interest and
importance. When $X$ is an affine space it was defined in
\cite{part1} (in even more specific situation of valuations
polynomial with respect to translations it was defined in
\cite{alesker-poly}). For a general manifold $X$ the
multiplicative structure was defined in \cite{part3}. The
construction in \cite{part3} uses the affine case \cite{part1} and
additional tools from the geometric measure theory.

It was shown in \cite{part3} that the product $V^\infty(X)\times
V^\infty(X)\to V^\infty(X)$ is a continuous map, and $V^\infty(X)$
becomes a commutative associative algebra with the unit (which is
the Euler characteristic).  The goal of this article is to study
further properties of the multiplicative structure and apply one
of them (which we call the Selfduality property) to introduce a
new class of generalized valuations.

In \cite{part2} a filtration of $V^\infty(X)$
\begin{eqnarray}\label{fltr}
V^\infty(X)=W_0\supset W_1\supset\dots\supset W_n
\end{eqnarray}
by closed subspaces was introduced. The first main result of this
article (Theorem \ref{comp1}) says that this filtration is
compatible with the product, namely $W_i\cdot W_j\subset W_{i+j}$
(where $W_k=0$ for $k>n$).

In \cite{part2} the author has introduced a continuous involution
$\sigma\colon V^\infty(X)\to V^\infty(X)$ called the Euler-Verdier
involution. The second main result of this article says that
$\sigma$ is an algebra automorphism (Theorem \ref{euler-verdier}).

Let us denote by $V^\infty_c(X)$ the space of compactly supported
smooth valuations. Next we introduce in this article the
integration functional $\int\colon V^\infty_c(X)\to \CC$. Slightly
oversimplifying, it is defined by $[\phi\mapsto \phi(X)]$. The
third main result is as follows.
\begin{theorem}\label{selfduality-property}
Consider the bilinear form
$$V^\infty(X)\times V^\infty_c(X)\to \CC$$
given by $(\phi,\psi)\mapsto \int \phi\cdot \psi$.

This bilinear form is a perfect pairing. More precisely the
induced map
$$V^\infty(X)\to (V^\infty_c(X))^*$$
is injective and has a dense image with respect to the weak
topology on $(V^\infty_c(X))^*$.
\end{theorem}
This is Theorem \ref{s1} in the paper. Its proof uses the
Irreducibility Theorem from \cite{alesker-gafa} in full
generality. Roughly Theorem \ref{selfduality-property} can be
interpreted as a selfduality property of the space of valuations
(at least when the manifold $X$ is compact).

Let us denote $V^{-\infty}(X):=(V^\infty_c(X))^*$. We call
$V^{-\infty}(X)$ the space of generalized valuations. We show
(Proposition \ref{2-2}) that $V^{-\infty}(X)$ has a canonical
structure of $V^\infty(X)$-module.

In \cite{part2} it was shown that the assignment to any open
subset $U\subset X$
$$U\mapsto V^\infty(U)$$
with the natural restriction maps is a sheaf denoted by
$\cv^\infty_X$. Here we show that
$$U\mapsto V^{-\infty}(U)$$
with the natural restriction maps is also a sheaf which we denote
by $\cv^{-\infty}_X$. Moreover $\cv^{-\infty}_X$ is a sheaf of
$\cv^\infty_X$-modules (Proposition \ref{2-9}).

Remind that by \cite{part2} the last term $W_n$ of the filtration
(\ref{fltr}) coincides with the space $C^\infty(X,|\ome_X|)$ of
smooth densities on $X$ (where $|\ome_X|$ denotes the line bundle
of densities on $X$), and $V^\infty(X)/W_1$ is canonically
isomorphic to the space of smooth functions $C^\infty(X)$. In
Subsection \ref{filtration} of this article we extend the
filtration $\{W_\bullet\}$ to generalized valuations by taking the
closure of $W_i$ in the weak topology on $V^{-\infty}(X)$:
$$V^{-\infty}(X)=W_0(V^{-\infty}(X))\supset
W_1(V^{-\infty}(X))\supset\dots\supset W_n(V^{-\infty}(X)).$$ We
show that $W_n(V^{-\infty}(X))$ is equal to the space
$C^{-\infty}(X,|\ome_X|)$ of generalized densities on $X$
(Proposition \ref{f6}). It is also shown that
$V^{-\infty}(X)/W_1(V^{-\infty}(X))$ is canonically isomorphic to
the space $C^{-\infty}(X)$ of generalized valuations on $X$
(Proposition \ref{f7}).


The Euler-Verdier involution is extended by continuity in the weak
topology to the space of generalized valuations (Subsection
\ref{euler-verdier-generalized}). Also the integration functional
extends (uniquely) by continuity in an appropriate topology to
generalized valuations with compact support (Subsection
\ref{euler-verdier-generalized}).

In Section \ref{real-analytic} we consider valuations on a real
analytic manifold $X$. On such a manifold one has the algebra of
constructible functions $\cf(X)$ which is a quite well known
object (see \cite{kashiwara-schapira}, Ch. 9). We construct a
canonical imbedding of the space $\cf(X)$ to the space of
generalized valuations $\vmi$ as a dense subspace. It turns out to
be possible to interpret some properties of valuations in more
familiar terms of constructible functions. Thus we show that the
canonical filtration on $\vmi$ induces on $\cf(X)$ the filtration
by codimension of support (Proposition \ref{ff3}). The restriction
of the integration functional to the space of compactly supported
constructible functions coincides with the well known functional
of integration with respect to the Euler characteristic
(Proposition \ref{gg1}). The restriction of the Euler-Verdier
involution on $\vmi$ to $\cf(X)$ coincides (up to a sign) with the
well known Verdier duality operator (Proposition \ref{vd2}).

{\bf Acknowledgements.} I express my gratitude J. Bernstein for
numerous stimulating discussions. I thank V.D. Milman for his
attention to this work. I thank A. Bernig for sharing with me the
recent preprint \cite{bernig-brocker}, J. Fu for very helpful
explanations on the geometric measure theory, P.D. Milman for
useful correspondences regarding subanalytic sets, and P. Schapira
for useful discussions on constructible sheaves and functions.

\section{Background.}\label{background-section} In this section we
fix some notation and remind various known facts. This section
does not contain new results.

In Subsection \ref{notation} we fix some notation and remind the
notions of normal and characteristic cycles of {\itshape convex}
sets. In Subsection \ref{suban} we review basic facts on
subanalytic sets. Subsection \ref{chc} collects facts on normal
and characteristic cycles. In Subsection \ref{some-valuations} we
review some notions on valuations on manifolds following mostly
\cite{part1}, \cite{part2}, \cite{part3}. Subsection
\ref{filtration-section} is also on valuations, and it reviews the
canonical filtration on valuations following \cite{part2}.
\subsection{Notation.}\label{notation}
Let $V$ be a finite dimensional real vector space.

$\bullet$ Let $\ck(V)$ denote the family of convex compact subsets
of $V$.

$\bullet$ Let $\RR_{\geq 0}$ (resp. $\RR_{>0}$) denote the set of
non-negative (resp. positive) real numbers.

$\bullet$ For a manifold $X$ let us denote by $|\ome_X|$ the line
bundle of densities over $X$.

$\bullet$ For a smooth manifold $X$ let $\cp(X)$ denote the family
of all simple subpolyhedra of $X$. (Namely $P\in \cp(X)$ iff $P$
is a compact subset of $X$ locally diffeomorphic to $\RR^k\times
\RR^{n-k}_{\geq 0}$ for some $0\leq k\leq n$. For a precise
definition see \cite{part2}, Subsection 2.1.)

$\bullet$ We denote by $\PP_+(V)$ the {\itshape oriented
projectivization} of $V$. Namely $\PP_+(V)$ is the manifold of
oriented lines in $V$ passing through the origin.

$\bullet$ For a vector bundle $E$ over a manifold $X$ let us
denote by $\PP_+(E)$ the bundle over $X$ whose fiber over any
point $x\in X$ is equal to $\PP_+(E_x)$ (where $E_x$ denotes the
fiber of $E$ over $x$).

$\bullet$ For a convex compact set $A\in \ck(V)$ let us denote by
$h_A$ the {\itshape supporting functional} of $A$, $h_A\colon
V^*\to \RR$. It is defined by
$$h_A(y):=sup\{y(x)| x\in A\}.$$

$\bullet$ Let $L$ denote the (real) line bundle over $\PP_+(V^*)$
such that its fiber over an oriented line $l\in \PP_+(V^*)$ is
equal to the dual line $l^*$.

$\bullet$ For a smooth vector bundle $E$ over a manifold $X$ and
$k$ being a non-negative integer or infinity, let us denote by
$C^k(X,E)$ the space of $C^k$-smooth sections of $E$. We denote by
$C^k_c(X,E)$ the space of $C^k$-smooth sections with compact
support. Let us denote by $C^{-\infty}(X,E)$ the space of
generalized sections of $E$ which is equal by definition to the
dual space $(C^\infty_c(X,E^*\otimes |\ome_X|))^*$. We have the
canonical imbedding $C^k(X,E)\hookrightarrow C^{-\infty}(X,E)$
(see e.g. \cite{guillemin-sternberg}, Ch. VI \S 1).


Let $K\in \ck(V)$. Let $x\in K$.
\begin{definition}\label{cc-1-1}
A tangent cone to $K$ at $x$ is the set denoted by $T_xK$ which is
equal to the closure of the set $\{y\in V|\exists \eps>0\, \,
x+\eps y\in K\}$.
\end{definition}
It is easy to see that $T_xK$ is a closed convex cone.
\begin{definition}\label{cc-1-2}
A normal cone to $K$ at $x$ is the set
$$(T_xK)^o:=\{y\in V^*| \,\, y(x)\geq 0 \forall x\in T_xK\}.$$
\end{definition}
Thus $(T_xK)^o$ is also a closed convex cone.
\begin{definition}\label{cc-1-3}
Let $K\in \ck(V)$. The {\itshape characteristic cycle} of K is the
set $$CC(K):=\cup_{x\in K}(T_xK)^o.$$
\end{definition}
It is easy to see that $CC(K)$ is a closed $n$-dimensional subset
of $T^*V=V\times V^*$ invariant with respect to the multiplication
by non-negative numbers acting on the second factor.
\def\ucc{\underline{CC}}
\def\tcc{\tilde{CC}}
For some references on the characteristic and normal cycles of
various sets see Remark \ref{remark-charact} below.

\subsection{Subanalytic sets.}\label{suban}
In this subsection we review some basic facts from the theory of
subanalytic sets of Hironaka. For more information see
\cite{hironaka1}, \cite{hironaka2}, \cite{hardt1}, \cite{hardt2},
\cite{bierstone-milman}, \cite{tamm}, and \S 8.2 of
\cite{kashiwara-schapira}. Let $X$ be a real analytic manifold.
\begin{definition}\label{sg1}
Let $Z$ be a subset of the manifold $X$. $Z$ is called {\itshape
subanalytic at a point $x\in X$} if there exists an open
neighborhood $U$ of $x$, compact real analytic manifolds $Y_j^i,\,
i=1,2,\, j=1,\dots,N$, and real analytic maps
$$f_j^i\colon Y_j^i\to X$$
such that
$$Z\cap U=U\cap \cup_{j=1}^N(f_j^1(Y_j^1)\backslash
f_j^2(Y_j^2)).$$

$Z$ is called {\itshape subanalytic in $X$} if $Z$ is subanalytic
at every point of $X$.
\end{definition}

\begin{proposition}\label{sg2}
(i) Let $Z$ be a subanalytic subset of the manifold $X$. Then the
closure and the interior of $Z$ are subanalytic subsets.

(ii) The connected components of a subanalytic set are locally
finite and subanalytic.

(iii) Let $Z_1$ and $Z_2$ be subanalytic subsets of the manifold
$X$. Then $Z_1\cup Z_2,\, Z_1\cap Z_2$, and $Z_1\backslash Z_2$
are subanalytic.
\end{proposition}

\begin{definition}\label{sg3}
Let $Z$ be a subanalytic subset of the manifold $X$. A point $x\in
Z$ is called {\itshape regular} if there exists an open
neighborhood $U$ of $x$ in $X$ such that $U\cap Z$ is a
submanifold of $X$.
\end{definition}
The set of regular points is denoted by $Z_{reg}$. Define the set
of {\itshape singular} points by $Z_{sing}:=Z\backslash Z_{reg}$.

\begin{proposition}\label{sg4}
The sets $Z_{reg}$ and $Z_{sing}$ are subanalytic, and $Z\subset
\bar Z_{reg}$.
\end{proposition}

If $x\in Z_{reg}$ then the dimension of $Z$ at $x$ is well
defined; it is denoted by $\dim_xZ$. Define $$\dim Z:=\sup_{x\in
Z_{reg}}\dim_x(Z).$$ Clearly $\dim Z\leq \dim X$.


\begin{proposition}\label{sg5}
Let $Z\subset X$ be a subanalytic subset. Then

(i) $\dim(Z\backslash Z_{reg})<\dim Z$;

(ii) $\dim (\bar Z\backslash Z)<\dim Z$.
\end{proposition}


\begin{definition}[\cite{kashiwara-schapira}, \S9.7]\label{construct-function} An integer valued function $f\colon
X\to \ZZ$ is called constructible if

1) for any $m\in \ZZ$ the set $f^{-1}(m)$ is subanalytic;

2) the family of sets $\{f^{-1}(m)\}_{m\in \ZZ}$ is locally
finite.
\end{definition}

Clearly the set of constructible $\ZZ$-valued functions is a ring
with pointwise multiplication. As in \cite{kashiwara-schapira} we
denote this ring by $CF(X)$. Define
\begin{eqnarray}\label{constr-def}
\cf:=CF(X)\otimes_\ZZ \CC.
\end{eqnarray}
Thus $\cf$ is a subalgebra of the $\CC$-algebra of complex valued
functions on $X$. In the rest of the article the elements of $\cf$
will be called {\itshape constructible functions}.

\def\cfx{\cf(X)}
 Let $\cf_c(X)$ denote the subspace of $\cf(X)$ of
{\itshape compactly supported} constructible functions. Clearly
$\cf_c(X)$ is a subalgebra of $\cf(X)$ (without unit if $X$ is
non-compact).

For a subset $P\subset X$ let us denote by $\One_P$ the indicator
function of $P$, namely
\begin{eqnarray*}
\One_P(x)=\left\{\begin{array}{ccc}
                       1&\mbox{ if }&x\in P\\
                       0&\mbox{ if }&x\not\in P.
                  \end{array}\right.
\end{eqnarray*}
\begin{proposition}\label{sg6}
(i) Any function $f\in \cf(X)$ can be presented locally as finite
linear combination of functions of the form $\One_Q$ where $Q$ is
a closed subanalytic subset.

(ii) Any function $f\in \cf_c(X)$ can be presented as finite
linear combination of functions of the form $\One_Q$ where $Q$ is
a compact subanalytic subset.
\end{proposition}
{\bf Proof.} Both statements are proved similarly. Let prove say
the second one. Let $f\in \cf_c(X)$. We prove the statement by the
induction on $\dim (\supp f)$ (note that $\supp f$ is a
subanalytic subset). If $\dim (\supp f)=0$ then there is nothing
to prove. Let us assume that we have proven the results for all
constructible functions with the dimension of support strictly
less than $k$. Let us prove it for $k$. Clearly $f$ is a finite
linear combination of functions of the form $\One_Q$ where $Q$ is
relatively compact subanalytic subset with $\dim Q\leq k$. But
$$\One_Q=\One_{\bar Q}-\One_{\bar Q\backslash Q}.$$
By Proposition \ref{sg2} the set $\bar Q\backslash Q$ is
subanalytic, and by Proposition \ref{sg5}(ii) $\dim(\bar
Q\backslash Q)<k$. The induction assumption implies the result.
\qed

\subsection{Characteristic and normal cycles.}\label{chc}
In Subsection \ref{notation} we have reminded the notion of
characteristic cycle of  convex compact sets. In this subsection
we remind the notion of characteristic cycle and very similar
notion of normal cycles of sets either from the class $\cp(X)$ on
a smooth manifold $X$, or the class of subanalytic subsets of a
real analytic manifold $X$ (in fact in the real analytic situation
these notions will be discussed more generally for constructible
functions on $X$ following \cite{kashiwara-schapira}). The notions
of characteristic and normal cycles of various classes of sets
coincide on the pairwise intersections of these classes.
\begin{remark}\label{remark-charact}
The notion of the characteristic cycle is not new. First an almost
equivalent notion of normal cycle (see below) was introduced by
Wintgen \cite{wintgen}, and then studied further by Z\"ahle
\cite{zahle87} by the tools of geometric measure theory.
Characteristic cycles of subanalytic sets of real analytic
manifolds were introduced by J. Fu \cite{fu-94} using the tools of
geometric measure theory and independently by Kashiwara (see
\cite{kashiwara-schapira}, Chapter 9) using the tools of the sheaf
theory. J. Fu's article \cite{fu-94} develops a more general
approach to define the normal cycle for more general sets than
subanalytic or convex ones (see Theorem 3.2 in \cite{fu-94}).
Applications of the method of normal cycles to integral geometry
can be found in \cite{fu-kinematic}.
\end{remark}
For simplicity of the exposition, in the rest of this subsection
we will assume that the manifold $X$ is oriented. Then
characteristic (resp. normal) cycle is a cycle in $T^*X$ (resp.
$\PP_+(T^*X)$. Nevertheless the characteristic and normal cycles
can be defined on non-oriented (even non-orientable) manifolds;
then they are cycles taking values is the local system $p^*o$
where $o$ is the orientation bundle over $X$ and $p\colon T^*X\to
X$ is the canonical projection. We refer to
\cite{kashiwara-schapira}, \S9.3, for the details on that. Though
in our applications to valuations of these notions we will need
the general case of not necessarily orientable manifolds, we will
ignore here this subtlety. Thus here we discuss the notions of
characteristic and normal cycles for oriented manifolds, but apply
it below for general manifolds.

Let us assume first that $X$ is a smooth oriented manifold. Set
$n=\dim X$. Let $P\in \cp(X)$. For any point $x\in P$ let us
define the {\itshape tangent cone} to $P$ at $x$, denoted by
$T_xP$, the set
$$T_xP:=\{\xi\in T_xX| \mbox{ there exists a } C^1-\mbox{map }
\gamma\colon [0,1]\to P \mbox{ such that }\gamma(0)=x\mbox{ and }
\gamma'(0)=\xi\}.$$ It is easy to see that $T_xP$ coincides with
the usual tangent space if $x$ is an interior point of $P$. In
general $T_xP$ is a closed convex polyhedral cone in $T_xX$.
Define
\begin{eqnarray}\label{j2}
CC(P):=\cup_{x\in P}(T_xP)^o
\end{eqnarray}
where for a convex cone $C$ in a linear space $W$ one denotes
$C^o$ its dual cone in $W^*$:
$$C^o:=\{y\in W^*|\, y(x)\geq 0\, \mbox{ for any } x\in C\}.$$
Clearly $CC(P)$ is invariant under the group $\RR_{>0}$ of
positive real numbers acting on the cotangent bundle $T^*X$ by
multiplication along the fibers. It is easy to see that $CC(P)$ is
an $n$-dimensional Lagrangian submanifold of $T^*X$ with
singularities. A choice of orientation on $X$ induces an
orientation on $CC(P)$. Then $CC(P)$ becomes a cycle, i.e. $\pt
(CC(P))=0$.

Let us assume now that $X$ is a {\itshape real analytic} manifold.
Again we assume that $X$ is oriented. Let $CF(X)$ be the ring of
integer valued constructible functions as in Definition
\ref{construct-function}, and let $\cf$ denote the algebra of
(complex valued) constructible functions as in (\ref{constr-def}).

In \cite{kashiwara-schapira}, \S9.7, there was constructed a group
homomorphism, also called characteristic cycle,
$$CC\colon CF(X)\to \cl(X)$$
where $\cl(X)$ denotes the group of Lagrangian conic subanalytic
cycles (with values in $p^*o$ in the non-oriented case). For the
formal definitions we refer to \cite{kashiwara-schapira}, \S\S
 9.7, 9.2. Here we describe $\cl(X)$ in a somewhat unformal way
when $X$ is oriented. An arbitrary element $\lam\in \cl(X)$ is an
$n$-cycle on $T^*X$ (i.e. $\pt \lam =0$) which locally over $X$
can be written as a finite sum $\lam=\sum_jm_j[\Lam_j]$ where
$m_j$ are integers, $\Lam_j$ are subanalytic oriented Lagrangian
locally closed submanifolds of $T^*X$ which are conic, i.e.
invariant under the action of the group of positive real numbers
$\RR_{>0}$ on $T^*X$, and $[\Lam_j]$ denotes the chain class of
$\Lam_j$.

Let us summarize some basic properties of $CC$ which will be used
later. First $CC$ commutes with restrictions of functions to open
subsets of $X$.

Let $P\subset X$ be a compact subanalytic subset. Assume in
addition that $P\in \cp(X)$. Then $CC(\One_P)$ coincides with the
characteristic cycle $CC(P)$ defined above in (\ref{j2}). Thus for
a subanalytic closed subset $Q$ we will also denote by $CC(Q)$ the
characteristic cycle $CC(\One_Q)$.

For a (locally closed) submanifold $S\subset X$ let us denote by
$T^*_SX$ the conormal bundle of $S$. If $S$ is subanalytic then
$T^*_SX$ is a subanalytic subset of $T^*X$ (Proposition 8.3.1 in
\cite{kashiwara-schapira}).

\begin{lemma}\label{lemma-charact}
Let $Q\subset X$ be a relatively compact subanalytic subset. Then
the closure $\bar Q$ can be presented as a finite union $\bar
Q=\cup_j Q_j$ of (locally closed) subanalytic submanifolds such
that
$$\supp (CC(\One_Q))\subset \cup_j T^*_{Q_j}X.$$
\end{lemma}
{\bf Proof.} Using induction in $\dim Q$ and Propositions
\ref{sg4}, \ref{sg5} we may replace $Q$ by $Q_{reg}$ and thus
assume that $Q$ is a (locally closed) submanifold of $X$.

Let us consider the subanalytic covering $X=Q\bigsqcup(X\backslash
Q)$. By Theorem 8.3.20 of \cite{kashiwara-schapira} there exists a
$\mu$-stratification $X=\bigsqcup_\beta X_\beta$ which is a
refinement of the above covering (for the definition of
$\mu$-stratification see Definition 8.3.19 of
\cite{kashiwara-schapira}).

Let us denote by $j\colon Q\to X$ the identity imbedding. Let
$\underline{\CC}_Q$ denote the constant sheaf on $Q$ (with complex
coefficients). Let $T_Q:=j_!\underline{\CC}_Q$ be the extension of
$\underline{\CC}_Q$ by zero. By the definition of the
characteristic cycle (\cite{kashiwara-schapira}, \S9.7)
\begin{eqnarray}\label{ks1}
CC(\One_Q)=CC(T_Q)
\end{eqnarray}
where in the right hand side stays the characteristic cycle of the
{\itshape sheaf} $T_Q$ (see \S 9.4 of \cite{kashiwara-schapira}).
Note that $T_Q$ is obviously constructible with respect to the
$\mu$-stratification $\{X_\beta\}$. It follows from the definition
of the characteristic cycle of a sheaf that
\begin{eqnarray}\label{ks2}
\supp CC(T_Q)\subset SS(T_Q)
\end{eqnarray}
where $SS(\cdot)$ denotes the singular support of a sheaf (see \S
5.1 of \cite{kashiwara-schapira}). Proposition 8.4.1 of
\cite{kashiwara-schapira} implies that $SS(T_Q)\subset
\bigsqcup_\beta T^*_{X_\beta}X$. But since $T_Q|_{X\backslash \bar
Q}=0$ one has
\begin{eqnarray}\label{ks3}
SS(T_Q)\subset \bigsqcup_{\beta\colon X_\beta\subset \bar
Q}T^*_{X_\beta}X.
\end{eqnarray}

Let us choose the covering $\bar Q=\cup_\alp Q_\alp$ where each
$Q_\alp$ is equal to one of the sets $X_\beta$ contained in $\bar
Q$. Thus (\ref{ks1})-(\ref{ks3}) imply
$$CC(Q)\subset \cup_\alp T^*_{Q_\alp}X.$$
Lemma is proved. \qed


Let us remind the definition of a normal cycle. We will treat all
the cases of subanalytic, convex, $\cp(X)$-sets, and constructible
functions simultaneously since in all these cases we already have
the notion of characteristic cycle.

Let $f$ be an element of one of these families. Let $CC(f)$ be its
characteristic cycle. Let us denote by $\underline{CC}(f)$ the
intersection of $CC(f)$ with the open subset of $T^*X$ obtained by
removing the zero section $\underline{0}$. Then
$\underline{CC}(f)$ is an $n$-cycle in $T^*X\backslash
\underline{0}$ invariant under the multiplication by positive real
numbers. Let $q\colon T^*X\backslash \underline{0}\to \PP_+(T^*X)$
denote the canonical quotient map. (Remind that $\PP_+(T^*X)$
denotes the bundle over $X$ whose fiber over a point $x\in X$ is
equal to the manifold of oriented lines in $T^*_xX$ passing
through the origin.)

It is easy to see that there exists unique $(n-1)$-cycle in
$\PP_+(T^*X)$ denoted by $\tilde{CC}(f)$ such that
$CC(f)=q^{-1}(\tilde{CC}(f))$. Consider the (antipodal) involution
$a\colon \PP_+(T^*X)\to \PP_+(T^*X)$ changing the orientation of
each line. Then by definition the normal cycle $N(f)$ is equal to
$a_*(\tilde{CC}(f))$. It is easy to see that if $CC(f)$ is a
subanalytic cycle then $N(f)$ is a subanalytic cycle, in
particular if $f$ is a constructible function then $N(f)$ is a
subanalytic cycle. Also it is known that $N(f)$ is a Legendrian
cycle when $\PP_+(T^*X)$ is equipped with the canonical contact
structure.

\subsection{Some valuation theory.}\label{some-valuations}
First let us remind some results from \cite{alesker-int}. Let $V$
be an $n$-dimensional real vector space. Let $\bar K=(K_1 , K_2
,\dots, K_s )$ be an $s$-tuple of compact convex subsets of $V $.
Let $r\in \NN\cup \{\infty\}$. For any $\mu \in C^{r}(V,\mov)$
consider the function $M_{\bar K} F \, :\RR^{s}_{+} \to \CC \,
,\mbox{where }\RR_{+}^{s}= \{(\lam_1,\dots,\lam_s) \,|
\,\lam_j\geq0 \} $, defined by    $$ ( M_{\bar K} \mu )
  (\lam _1,\dots,\lam _s)= \mu(\sum_{i=1}^{s} \lam _i K_i) .  $$
\begin{theorem}[\cite{alesker-int}]\label{minkowski}
(1) $M_{\bar K}\mu \in C^{r}(\RR^s_+)$ and $M_{\bar K}$ is a
continuous operator from $C^{r}(V,\mov)$ to $C^{r}(\RR^s_+)$.

(2) Assume that  a sequence $\mu^{(m)}$ converges to $\mu$ in
$C^{r}(V,\mov)$. Let $ K_{i}^{(m)}, \, K_{i},\,
i=1,\dots,s,\,m\in\NN $, be convex compact sets in $ V $, and for
every $ i=1,\dots ,s $ $K_{i}^{(m)} \to K_{i} $ in the Hausdorff
metric as $m\to \infty$. Then $ M_{\bar K^{(m)} } \mu^{(m)} \to
M_{\bar K } \mu $ in $C^{r}(\RR^s_+)$ as $m\to \infty$.
\end{theorem}

\begin{definition}
a) A function $\phi :{\cal K}(V) \to \CC$ is called a valuation if
for any $K_1, \, K_2 \in {\cal K}(V)$ such that their union is
also convex one has
$$\phi(K_1 \cup K_2)= \phi(K_1) +\phi(K_2) -\phi(K_1 \cap K_2).$$

b) A valuation $\phi$ is called continuous if it is continuous
with respect to the Hausdorff metric on ${\cal K}(V)$.
\end{definition}

For the classical theory of valuations we refer to the surveys
McMullen-Schneider \cite{mcmullen-schneider} and McMullen
\cite{mcmullen-survey}. For the general background from convexity
we refer to Schneider \cite{schneider-book}.

In \cite{part1} one has introduced a class $SV(V)$ of valuations
called {\itshape smooth valuations}. We refer to \cite{part1} for
an axiomatic definition. Here we only mention that $SV(V)$ is a
$\CC$-linear space (with the obvious operations) with a natural
Fr\'echet topology. In this article we will need a description of
$SV(V)$ which is Theorem \ref{onto} below.

Let us denote by ${}\!^ {\textbf{C}}L$ the (complex) line bundle
over $\PP_+(V^*)$ whose fiber over $l\in \PP_+(V^*)$ is equal to
$l^*\otimes_{\RR}\CC$ (where $l^*$ denotes the dual space to $l$).

Note that for any convex compact set $A\in \ck(V)$ the supporting
functional $h_A$ is a continuous section of ${}\!^ {\textbf{C}}L$,
i.e $h_A\in C(\PP_+(V^*),{}\!^ {\textbf{C}}L)$.
\begin{theorem}[\cite{part1}, Corollary 3.1.7]\label{onto}
There exists a continuous linear map
\begin{eqnarray*}
\ct\colon \oplus_{k=0}^nC^\infty(V\times
\PP_+(V^*)^k,|\ome_V|\boxtimes {}\!^ {\textbf{C}}L^{\boxtimes k})
\to SV(V)
\end{eqnarray*}
which is uniquely characterized by the following property: for any
$k=0,1,\dots,n$, any $\mu\in C^\infty(V,|\ome_V|)$, any strictly
convex compact sets $A_1,\dots,A_k$ with smooth boundaries, and
any $K\in \ck(V)$ one has
\begin{eqnarray*}
\ct(\mu\boxtimes h_{A_1}\boxtimes \dots\boxtimes h_{A_k})(K)=
\frac{\pt^k}{\pt\lam_1\dots \pt\lam_k}\big|_0
\mu(K+\sum_{i=1}^k\lam_iA_i)
\end{eqnarray*}
where $\lam_i\geq 0$ in the right hand side.

Moreover the map $\ct$ is an epimorphism.
\end{theorem}

In \cite{part2} one has introduced for any smooth manifold $X$ a
class of finitely additive measures on the family of simple
subpolyhedra $\cp(X)$. This class is denoted by $V^\infty(X)$. It
is a $\CC$-linear space (with the obvious operations). Then
$V^\infty(X)$ has a natural Fr\'echet topology. Moreover in the
case  of linear Fr\'echet space $V$ any element $\phi\in
V^\infty(V)$ being restricted to $\ck(V)\cap \cp(V)$ has a
(unique) extension by continuity in the Hausdorff metric to
$\ck(V)$, and this extension belongs to $SV(V)$. Thus one gets a
linear map
$$V^\infty(V)\to SV(V).$$
In \cite{part2}, Proposition 2.4.10, the following result was
proved.
\begin{proposition}\label{isomorphism}
The above constructed map $V^\infty(V)\to SV(V)$ is an isomorphism
of Fr\'echet spaces.
\end{proposition}

We will also need the following description of $V^\infty(X)$
obtained in \cite{part2} (based on some results on normal cycles
from Section 2 of \cite{part3}). Let us denote by $T^*X$ the
cotangent bundle of $X$. Let $p\colon T^*X\to X$ be the canonical
projection. Let $\Ome^n$ denote the vector bundle of $n$-forms
over $T^*X$. Let us denote by $o$ the orientation bundle over $X$.
Let us denote by $\tilde C^\infty(T^*X,\Ome^n\otimes p^*o)$ the
space of $C^\infty$-smooth sections of the bundle $\Ome^n\otimes
p^*o$ such that the restriction of the projection $p$ to the
support of such section is a proper map.
\begin{theorem}[\cite{part2}, Theorem 0.1.3]\label{normal cycle descr}
(i) Let $\ome \in\tilde C^\infty(T^*X,\Ome^n\otimes p^*o)$. The
functional $\cp(X)\to \CC$
\begin{eqnarray}\label{intfor}
P\mapsto \int_{CC(P)} \ome
\end{eqnarray}
is a smooth valuation.

(ii) Conversely, any smooth valuation $\phi\in V^\infty(X)$ has
the form (\ref{intfor}), i.e. there exists a form $\ome \in\tilde
C^\infty(T^*X,\Ome^n\otimes p^*o)$ such that $\phi(P)=\int_{CC(P)}
\ome$ for any $P\in \cp(X)$.
\end{theorem}
\begin{remark}
(1) The integration (\ref{intfor}) is well defined since a choice
of orientation of the manifold $X$ induces an orientation of
$CC(P)$.

(2) A presentation of a valuation $\phi$ in the form
(\ref{intfor}) is highly non-unique.
\end{remark}

Let us describe the multiplicative structure on $V^\infty(X)$
following \cite{part3}. It was shown in \cite{part2} that the
assignment to any open subset $U\subset X$ $$U\mapsto
V^\infty(U)$$ with the natural restriction maps is a sheaf. The
product on smooth valuations commutes with the restrictions to
open subsets. Hence it is enough to describe the product locally,
say under the assumption that $X$ is diffeomorphic to $\RR^n$. Let
us fix a diffeomorphism $X\tilde \to \RR^n$. Proposition
\ref{isomorphism} provides an isomorphism
$V^\infty(\RR^n)\tilde\to SV(\RR^n)$. In \cite{part1} the author
has described the product on $SV(\RR^n)$ which we will remind
below. The main point of \cite{part3} was to show that the
obtained product on $V^\infty(X)$ does not depend on the choices
of diffeomorphisms.

Thus it remains to describe the product on $SV(\RR^n)$ following
\cite{part1}. The product
$$SV(\RR^n)\times SV(\RR^n)\to SV(\RR^n)$$
is a continuous map which is uniquely defined by the
distributivity and the following property: let $\phi,\psi\in
SV(\RR^n)$ have the form
\begin{eqnarray}\label{b1}
\phi(K)=\frac{\pt^k}{\pt\lam_1\dots\pt\lam_k}\big|_0\mu(K+\sum_{i=1}^k\lam_iA_i),\\
\psi(K)=\frac{\pt^l}{\pt\mu_1\dots\pt\mu_l}\big|_0\mu(K+\sum_{j=1}^l\mu_jB_j)
\end{eqnarray}
where $0\leq k,l\leq n$; $\mu, \nu$ are smooth densities on
$\RR^n$; $A_1,\dots,A_k,B_1,\dots,B_l$ are strictly convex compact
sets with smooth boundaries, $K$ is an arbitrary convex compact
subset in $\RR^n$. Then
\begin{eqnarray}\label{b2}
(\phi\cdot\psi)(K)=\frac{\pt^{k+l}}{\pt\lam_1\dots\pt\lam_k\pt\mu_1\dots\pt\mu_l}\big|_0
(\mu\boxtimes\nu)\left(\Delta(K)+((\sum_{i=1}^k\lam_iA_i)\times
(\sum_{j=1}^l\mu_jB_j))\right)
\end{eqnarray}
where $\Delta\colon \RR^n \hookrightarrow \RR^n\times \RR^n$ is
the diagonal imbedding, $\mu\boxtimes \nu$ denotes the usual
exterior product of densities. Note that in (\ref{b1})-(\ref{b2})
the derivatives exist due to Theorem \ref{minkowski}.

Equipped with this product, the space $V^\infty(X)$ becomes a
commutative associative algebra with unit (the unit is the Euler
characteristic).

Let us describe the Euler-Verdier involution $\sigma$ on
$V^\infty(X)$ following \cite{part2}. Let $a\colon T^*X\to T^*X$
be the involution of multiplication by $-1$ in each fiber of the
projection $p\colon T^*X\to X$. It induces the involution
$$a^*\colon \tilde C^\infty(T^*X,\Ome^n\otimes p^*o)\to \tilde
C^\infty(T^*X,\Ome^n\otimes p^*o).$$ We have the following
proposition.
\begin{proposition}[\cite{part2}, Proposition 3.3.1]\label{ev}
The involution $(-1)^n a^*$ factorizes (uniquely) to the
involution of $V^\infty(X)$ which is denoted by $\sigma$. Moreover
$\sigma$ commutes with the restrictions to open subsets and thus
induces an involution of the sheaf $\cv^\infty_X$ which is also
denoted by $\sigma$.
\end{proposition}

\subsection{Filtration on valuations.}\label{filtration-section} In \cite{part2} we have
introduced on the space of smooth valuations $V^\infty(X)$ a
canonical finite filtration by closed subspaces:
\begin{eqnarray}\label{i1}
V^\infty(X)=W_0(X)\supset W_1(X)\supset \dots\supset W_n(X)
\end{eqnarray}
where $n=\dim X$. Let us remind some of the main properties of
this filtration.

\begin{proposition}[\cite{part2}, Proposition 3.1.2]\label{i2}
The assignment to each open subset $U\subset X$
$$U\mapsto W_i(U)$$
is a subsheaf of $\cv_X^\infty$. (This sheaf is denoted by
$\cw_i$.)
\end{proposition}
It turns out that the associated graded sheaf
$gr_\cw\cv^\infty_X:=\bigoplus_{i=0}^n\cw_i/\cw_{i+1}$ admits a
simple description in terms of translation invariant valuations.
To state it let us denote by $Val(TX)$ the (infinite dimensional)
vector bundle over $X$ such that its fiber over a point $x\in X$
is equal to the space $Val^{sm}(T_xX)$ of smooth translation
invariant valuations on the tangent space $T_xX$. By McMullen's
theorem \cite{mcmullen-euler} the space $Val^{sm}(T_xX)$ has
natural grading by the degree of homogeneity which must be an
integer between 0 and $n$. Thus $Val(TX)$ is a graded vector
bundle. Let us denote by $\underline{Val}(TX)$ the sheaf $U\mapsto
C^\infty(U,Val(TX))$ where the last space denotes the space of
infinitely smooth sections of $Val(TX)$ over $U$.
\begin{theorem}[\cite{part2}, Theorem 0.1.2 and Section 3]\label{i3}
There exists a canonical isomorphism of graded sheaves
$$gr_\cw\cv^\infty_X\simeq \underline{Val}(TX).$$
Moreover for any open subset $U\subset X$ the induced isomorphism
on global sections is isomorphism of linear topological spaces.
\end{theorem}
This theorem provides a description of smooth valuations since
translation invariant valuations are studied much better.
\begin{remark}\label{i4}
Interpreted appropriately, Theorem \ref{i3} says in particular
that the last term of the filtration $\cw_n$ is canonically
isomorphic to the sheaf of $C^\infty$-smooth measures (=densities)
on $X$, and the first quotient $\cv^\infty_X/\cw_1$ is canonically
isomorphic to the sheaf of $C^\infty$-smooth functions on $X$.
\end{remark}

The filtration $\{W_\bullet\}$ on valuations can be interpreted in
terms of Theorem \ref{normal cycle descr} as follows. First remind
the general construction of a filtration differential forms on a
total space of a bundle.

Let $X$ be a smooth manifold. Let $p:P\to X$ be a smooth bundle.
Let $\Omega^N(P)$ be the vector bundle over $P$ of $N$-forms. For
a vector space $R$ we denote by $Gr_N(R)$ the Grassmannian of
$N$-dimensional linear subspaces in $R$. Let us introduce a
filtration of $\Omega^N(P)$ by vector subbundles $W_i(P)$ as
follows. For every $y\in P$ set
\begin{multline*} (W_i(P))_y:=\{\omega\in \wedge^NT_y^*P \big|\,
\omega|_F\equiv 0 \mbox{ for all } F\in Gr_N(T_yP)\\
\mbox{ with } \dim(F\cap T_y(p^{-1}p(y)))>N-i\}.\end{multline*}

Clearly we have
$$\Omega^N(P)=W_0(P)\supset W_1(P)\supset \dots \supset W_N(P)\supset
W_{N+1}(P)=0.$$ Let us discuss this filtration in greater detail
following \cite{part1}.

Let us make some elementary observations from linear algebra.
\def\wel{W(L,E)}
Let $L$ be a finite dimensional vector space. Let $E\subset L$ be
a linear subspace. For a non-negative integer $i$ set
$$\wel_i:=\{\omega\in \wedge^NL^*\big|\, \omega|_F\equiv 0\mbox{ for
all } F\subset L \mbox{ with } \dim(F\cap E)>N-i\}.$$ Clearly
$$\wedge^NL^*=\wel_0\supset \wel_1\supset \dots \supset \wel_N\supset
\wel_{N+1}=0.$$
\begin{lemma}[\cite{part1}, Lemma 5.2.3] \label{val-19}
There exists canonical isomorphism of vector spaces
$$\wel_i/\wel_{i+1}=\wedge^{N-i}E^*\otimes \wedge^i(L/E)^*.$$
\end{lemma}

Let us apply this construction in the context of integration with
respect to the characteristic cycle. Let $X$ be a smooth manifold
of dimension $n$. Let $P:=T^*X$ be the cotangent bundle. Let
$p:P\to X$ be the canonical projection. Let us denote by $o$ the
orientation bundle on $X$.  The above construction gives a
filtration of $\Omega^n(P)$ by subbundles
$$\Omega^n(P)=W_0(\Omega^n(P))\supset\dots\supset W_n(\Omega^n(P)).$$
Twisting this filtration by $p^*o$ we get a filtration of
$\Omega^n(P)\otimes p^*o$ by subbundles denoted by
$W_i(\Omega^n(P)\otimes p^*o)$.

Let us denote by $\tilde C^\infty(P,W_i(\Omega^n\otimes p^*o))$
the space of infinitely smooth sections of the bundle
$W_i(\Omega^n\otimes p^*o)$ such that the restriction of the
projection $p$ to the support of these sections is proper. The we
have the following result.
\begin{theorem}[\cite{part2}, Proposition 3.1.9]
For any valuation $\phi\in W_i(X)$ there exists $\ome \in\tilde
C^\infty(T^*X,W_i(\Ome^n\otimes p^* o))$ such that for any $P\in
\cp(X)$
$$\phi(P)=\int_{CC(P)}\ome.$$
Conversely any such valuation belongs to $W_i(X)$.
\end{theorem}

\section{A technical
lemma.}\label{technical-lemma}\setcounter{subsection}{1}\setcounter{theorem}{0}
In this section we will prove a technical lemma which will be used
later on in this article.
\begin{lemma}\label{ccc4}
Let $i=0,1,\dots,n$. Let $\phi\in W_i(\vic)$. Then there exists a
compactly supported form $\ome\in
C^\infty_c(T^*X,W_i(\Ome^n\otimes p^*o))$ such that
$$\phi(P)=\int_{CC(P)}\ome \mbox{ for any } P\in \cp(X).$$
\end{lemma}
\begin{remark}\label{tech-remark}
A version of this lemma for smooth valuations without the
assumption on the compactness of support was proved in
\cite{part2}, Proposition 3.1.9; it will be used in the proof of
Lemma \ref{technical-lemma}.
\end{remark}
{\bf Proof} of Lemma \ref{technical-lemma}. As in \cite{part2}
consider the sheaves on $X$
\begin{eqnarray}\label{m1}
\cw_i(U)=W_i(V^\infty(U)),\\\label{m2} \cw_i'(U)=\tilde
C^\infty(T^*U,W_i(\Ome^n\otimes p^*o))
\end{eqnarray}
for any open subset $U\subset X$; in equality (\ref{m2}) the
symbol $\tilde C^\infty$ denotes the space of infinitely smooth
sections of a vector bundle over $T^*U$ such that the restriction
of the canonical projection $p\colon T^*U\to U$ to the support of
such sections is proper. We have the obvious inclusions:
\begin{eqnarray*}
\cw_n'\subset \cw_{n-1}'\subset \dots\subset \cw'_0;\\
\cw_n\subset \cw_{n-1}\subset \dots\subset \cw_0=\cv^\infty_X.
\end{eqnarray*}

The integration over the the characteristic cycle gives a morphism
of sheaves
\begin{eqnarray}\label{mm1}
T_i\colon \cw'_i\to \cw_i.
\end{eqnarray}
By Proposition 3.1.9 of \cite{part2} $T_i$ is an epimorphism of
sheaves. Clearly the restriction of $T_i$ to $\cw'_{i+1}$ is equal
to $T_{i+1}$. Define the sheaves
\begin{eqnarray}\label{mm2}
\ck_i:=Ker T_i.
\end{eqnarray}
We obviously have
\begin{eqnarray*}
\ck_n\subset \ck_{n-1}\subset \dots\subset \ck_0\subset \cw'_0.
\end{eqnarray*}
Let us consider the associated graded sheaves
\begin{eqnarray}\label{mm3}
\cf:=\oplus_{i=o}^n\cw_i/\cw_{i+1},\\\label{mm4}
\cf':=\oplus_{i=0}^n\cw'_i/\cw'_{i+1}.
\end{eqnarray}
The epimorphism $T_0\colon \cw'_0\twoheadrightarrow \cw_0$ induces
the epimorphism
\begin{eqnarray}\label{mm4.5}
T\colon \cf'\twoheadrightarrow \cf.
\end{eqnarray}
Define $\ct:=Ker T$. Clearly
\begin{eqnarray}\label{mm5}
\ct=\oplus_{i=0}^n\ck_i/\ck_{i+1}.
\end{eqnarray}
Let us denote by $\co_X$ the sheaf of $C^\infty$-smooth functions
on $X$. It was shown in \cite{part2} (see the proof of Proposition
3.1.9) that $\ct$ is naturally isomorphic to a sheaf of
$\co_X$-modules. Hence by Section 3.7 of Ch. II in Godement's book
\cite{godement} one has
\begin{eqnarray*}
H^j_c(X,\ct)=0 \mbox{ for } j>0.
\end{eqnarray*}
Hence $H_c^j(X,\ck_i/\ck_{i+1})=0$ for $j>0,\, i=0,1,\dots,n$. By
the long exact sequence we deduce
\begin{eqnarray}\label{mm6}
H^j_c(X,\ck_i)=0 \mbox{ for } j>0,\, i=0,1,\dots,n.
\end{eqnarray}
For the short exact sequence of sheaves
\begin{eqnarray*}
0\to\ck_i\to \cw_i'\overset{T_i}{\to}\cw_i\to 0
\end{eqnarray*}
consider the beginning of the long exact sequence in cohomology
with compact support
\begin{eqnarray}\label{mm7}
H^0_c(X,\cw'_i)\to H^0_c(X,\cw_i)\to H^1_c(X,\ck_i).
\end{eqnarray}
But the last space in (\ref{mm7}) vanishes due to (\ref{mm6}).
Hence the map $H^0_c(X,\cw'_i)\to H^0_c(X,\cw_i)$ is surjective.
But
\begin{eqnarray*}
H^0_c(X,\cw_i)=W_i(\vic);\\
H^0_c(X,\cw'_i)=C^\infty_0(T^*X,\Ome^n\otimes p^*o).
\end{eqnarray*}
Thus lemma is proved. \qed

\section{Compatibility of the filtration with the
product.}\label{compatibility}
\setcounter{subsection}{1}\setcounter{theorem}{0}\setcounter{equation}{0}
The main results of this section are Theorems \ref{comp1} and
\ref{comp4} below.

Remind that in Subsection \ref{filtration-section} we have
discussed the canonical filtration by closed subspaces
$$V^\infty(X)=W_0(X)\supset W_1(X)\supset\dots\supset W_n(X).$$
It will be convenient to extend this filtration infinitely by
putting
$$W_i(X)=0 \mbox{ for } i>n.$$

\begin{theorem}\label{comp1}
For any $i,j\geq 0$ one has
$$W_i(X)\cdot W_j(X)\subset W_{i+j}(X).$$
\end{theorem}
{\bf Proof.} By Corollary 4.1.4 of \cite{part3} $\cv^\infty_X$ is
a sheaf of algebras, i.e. the product commutes with the
restriction to open subsets. Hence we may assume that $X$ is
diffeomorphic to $\RR^n$. Let us fix a diffeomorphism $X\tilde\to
\RR^n$. Let us consider the induced isomorphism of linear
topological spaces
$$V^\infty(X)\tilde\to SV(\RR^n)$$
from Proposition \ref{isomorphism}. By Proposition 3.1.3 of
\cite{part2} the subspace $W_i(X)$ is isomorphic under this
isomorphism to a closed subspace of $SV(\RR^n)$. Let us denote
this subspace by $\hat W_i$; it was explicitly defined in Section
3 of \cite{part1} in slightly different notation. Moreover by
Theorem 4.1.2(4) of \cite{part1}
\begin{eqnarray}\label{comp2}
\hat W_i\cdot \hat W_j\subset \hat W_{i+j}.
\end{eqnarray}
Hence our result follows from (\ref{comp2}) and the construction
of the product on $V^\infty(X)$ described in Subsection
\ref{some-valuations}. \qed

Remind that by Theorem \ref{i3} there exists a canonical
isomorphism of graded linear topological spaces
\begin{eqnarray}\label{comp3}
gr_WV^\infty(X):=\oplus_{i=0}^nW_i(X)/W_{i+1}(X)\tilde\to
C^\infty(X,Val(TX))
\end{eqnarray}
where the vector bundle $Val(TX)$ over $X$ was defined in
Subsection \ref{filtration-section}.

Observe that $gr_W\vi$ is a graded algebra with the product
induced from $\vi$. Note also that $C^\infty(X,Val(TX))$ is also a
graded algebra with the product defined pointwise. Namely if
$f,g\in C^\infty(X,Val(TX))$ then for any point $x\in X$
$$(f\cdot g)(x)=f(x)\cdot g(x)\in Val^{sm}(T_xX).$$
We are going to prove
\begin{theorem}\label{comp4}
The isomorphism (\ref{comp3}) is an isomorphism of algebras.
\end{theorem}
{\bf Proof.} As in the proof of Theorem \ref{comp1}, the statement
is local. Thus we may assume that $X$ is diffeomorphic to $\RR^n$.
Now the result follows from the construction of the product
described in Subsection \ref{some-valuations} and Theorem 4.1.3 of
\cite{part1} where the corresponding statement was proved for
valuations on convex subsets of a linear space. \qed

\section{The automorphism property of the Euler-Verdier
involution.}\label{euler-verdier-section}
\setcounter{subsection}{1}\setcounter{theorem}{0}\setcounter{equation}{0}
The main result of this section is Theorem \ref{euler-verdier}.
\begin{lemma}\label{l0.1}
Let $\phi\in V^\infty(X)$. Let $P\in \cp(X)$. Then
$$(\sigma\phi)(P)=(-1)^{\dim P} (\phi(P)-\phi(\pt P)).$$
\end{lemma}
{\bf Proof.} Equality (15) in \cite{part2} says that for any
$\ome\in \tilde C^\infty(T^*X,\Ome^n\otimes p^*o)$ and any $P\in
\cp(X)$ one has
\begin{eqnarray}\label{ur}
\int_{CC(P)}a^* \ome=(-1)^{n-\dim P}\left(\int_{CC(P)}\ome
-\int_{CC(\pt P)}\ome\right)
\end{eqnarray}
where $\pt P:=P\backslash int P$ and $int P$ denotes the relative
interior of $P$. The result follows immediately from Proposition
\ref{ev} and (\ref{ur}). \qed

From Lemma \ref{l0.1} we immediately deduce that the Euler-Verdier
involution commutes with restriction to submanifolds. More
precisely we have the following lemma.
\begin{lemma}\label{l0.2}
Let $Y$ be a smooth submanifold of a manifold $X$. Let $\phi\in
V^\infty(X)$. Then
$$(\sigma\phi)|_Y=\sigma(\phi|_Y).$$
\end{lemma}

\begin{lemma}\label{l1}
Let $\phi\in V^\infty(\RR^n)$ be a smooth valuation such that for
any $K\in \ck(\RR^n)$ one has
\begin{eqnarray*}
\phi(K)=\frac{\pt^k}{\pt\lam_1\dots \pt\lam_k}\big
|_0\mu(K+\sum_{i=1}^k\lam_iA_i)
\end{eqnarray*}
where $\mu$ is a smooth density on $\RR^n$, and $A_1,\dots,A_k$
are strictly convex compact subsets with smooth boundaries and
containing the origin in the interior. Then
\begin{eqnarray}\label{e1}
(\sigma\phi)(K)=(-1)^{n-k}\frac{\pt^k}{\pt\lam_1\dots
\pt\lam_k}\big |_0\mu(K+\sum_{i=1}^k\lam_i(-A_i)).
\end{eqnarray}
\end{lemma}
{\bf Proof.} For $k=0$ the lemma is obvious. Let us assume that
$k>0$. It is enough to prove (\ref{e1}) under the assumption that
$K$ has non-empty interior and strictly convex smooth boundary.
For any $\lam_1,\dots,\lam_k>0$ the map
\begin{eqnarray*}
\Xi_{\lam_1\dots\lam_k}\colon V\times \PP_+(V^*)\times (0,1]\to V
\end{eqnarray*}
given by $(p,n,t)\mapsto p+t\sum_{i=1}^k\lam_i\nabla h_{A_i}(n)$
induces a homeomorphism of $N(K)\times (0,1]$ onto its image
$(K+\sum_{i=1}^k\lam_iA_i)\backslash K$ (this is well known; see
e.g. Proposition 3.1.2 of \cite{part3} where this statement is
proved under some more general assumptions). Hence
\begin{eqnarray*}
\phi(K)=\frac{\pt ^k}{\pt\lam_1\dots
\pt\lam_k}\big|_0\int_{N(K)\times[0,1]}
\Xi_{\lam_1\dots\lam_k}^*\mu.
\end{eqnarray*}

Let us denote by $\tilde a\colon \PP_+(V^*)\to \PP_+(V^*)$ the
involution of changing an orientation of a line. Then
\begin{eqnarray}\label{e2}
(\sigma\phi)(K)=(-1)^n\frac{\pt^k}{\pt\lam_1\dots\pt\lam_k}\big|_0
\tilde a^*(\Xi_{\lam_1\dots\lam_k}^*\mu)=(-1)^n
\frac{\pt^k}{\pt\lam_1\dots\pt\lam_k}\big|_0
((\Xi_{\lam_1\dots\lam_k}\circ \tilde a)^*\mu).
\end{eqnarray}
Observe that
\begin{eqnarray*}
(\Xi_{\lam_1\dots\lam_k}\circ
a)(p,n,t)=p+t\sum_{i=1}^k\lam_i(\nabla h_{A_i})(-n).
\end{eqnarray*}
But $h_{-A}(n)=h_A(-n)$. Hence
\begin{eqnarray}\label{e3}
(\Xi_{\lam_1\dots\lam_k}\circ
a)(p,n,t)=p-t\sum_{i=1}^k\lam_i(\nabla h_{-A_i})(n).
\end{eqnarray}
Note that
\begin{eqnarray}\label{e4}
\frac{\pt^k}{\pt\lam_1\dots\pt\lam_k}\big|_0\mu(K+\sum_{i=1}^k\lam_i(-A_i))=
\frac{\pt^k}{\pt\lam_1\dots\pt\lam_k}\big|_0\int_{N(K)\times[0,1]}
\tilde\Xi_{\lam_1\dots\lam_k}\mu
\end{eqnarray}
where $\tilde\Xi_{\lam_1\dots\lam_k}\colon
V\times\PP_+(V^*)\times[0,1]\to V$ is defined by
$\tilde\Xi_{\lam_1\dots\lam_k}(p,n,t)=p+t\sum_{i=1}^k \lam_i\nabla
h_{-A_i}(n)$. Now Lemma \ref{l1} follows from
(\ref{e2}),(\ref{e3}),(\ref{e4}). \qed

\begin{theorem}\label{euler-verdier}
The Euler-Verdier involution $\sigma\colon \cv^\infty_X\to
\cv^\infty_X$ is an algebra automorphism. Moreover it preserves
the filtration $\cw_\bullet$, namely $\sigma(\cw_i)=\cw_i$ for any
$i=0,\dots,n$.
\end{theorem}
{\bf Proof.} The second part of the theorem was proved in
\cite{part2}. Thus it remains to show that $\sigma$ is an algebra
automorphism. The statement is local thus we may and will assume
that $X=\RR^n$. Let $\phi,\psi\in V^\infty(\RR^n)$. We may assume
that for any $K\in \ck(\RR^n)$
\begin{eqnarray*}
\phi(K)=\frac{d^k}{d\eps^k}\big|_0\mu(K+\eps A),\,\,
\psi(K)=\frac{d^l}{d\delta^l}\big|_0\nu(K+\delta B)
\end{eqnarray*}
where $\mu,\nu$ are smooth densities on $\RR^n$, and $A,B$ are
strictly convex compact subsets with smooth boundaries and
containing the origin the interior. Then
\begin{eqnarray*}
(\phi\cdot \psi)(K)=\frac{\pt^{k+l}}{\pt^k\eps\cdot
\pt^l\delta}\big|_0(\mu\boxtimes \nu)(\Delta(K)+(\eps A,\delta B))
\end{eqnarray*}
where $\Delta\colon \RR^n \hookrightarrow \RR^n\times \RR^n$ is
the diagonal imbedding. By Lemma \ref{l1} one has
\begin{eqnarray*}
(\sigma \phi)(K)=(-1)^{n-k}\frac{d^k}{d\eps^k}\big|_0\mu(K+\eps
(-A)),\,\, (\sigma
\psi)(K)=(-1)^{n-l}\frac{d^l}{d\delta^l}\big|_0\nu(K+\delta
(-B))\\
(\sigma(\phi\cdot
\psi))(K)=(-1)^{2n-(k+l)}\frac{\pt^{k+l}}{\pt^k\eps\cdot
\pt^l\delta}\big|_0(\mu\boxtimes \nu)(\Delta(K)+(\eps (-A),\delta
(-B))).
\end{eqnarray*}
Hence we have
\begin{eqnarray*}
(\sigma\phi\cdot
\sigma\psi)(K)=(-1)^{2n-(k+l)}\frac{\pt^{k+l}}{\pt^k\eps\cdot
\pt^l\delta}\big|_0(\mu\boxtimes \nu)(\Delta(K)+(\eps (-A),\delta
(-B)))=(\sigma(\phi\cdot \psi))(K).
\end{eqnarray*}
\qed

\section{The integration functional on
valuations.}\label{integration-section} In Subsection
\ref{val-comp-supp} we describe canonical linear topology on the
space $\vic$ of compactly supported smooth valuations. In
Subsection \ref{integration-functional} we construct a canonical
continuous linear functional $\vic\to \CC$ called the {\itshape
integration functional}.
\subsection{Valuations with compact support.}\label{val-comp-supp} In this subsection we introduce
the space of valuations $V^\infty_c(X)$ with compact support and
establish some of the simplest properties of it.

Let $\phi\in V^\infty(X)$. We say that a point $x\in X$ does not
belong to the support of $\phi$ if there exists a neighborhood $U$
of $x$ such that $\phi|_U\equiv 0$. The set of all points which
does not belong to support of $\phi$ is an open subset of $X$. Its
complement is called the support of $\phi$ and is denoted by
$\supp\, \phi$. Thus $\supp\, \phi$ is a closed subset of $X$. The
following lemma is obvious.
\begin{lemma}\label{c1}
For any $\phi,\psi\in V^\infty(X)$
$$\supp \,(\phi\cdot \psi)\subset \supp\, \phi\cap \supp \,\psi.$$
\end{lemma}
The space of all valuations with compact support will be denoted
by $V^\infty_c(X)$. Also for any subset $S\subset X$ let us denote
$$V_S^\infty(X):=\{\phi\in V^\infty(X)\, |\, \supp \,\phi\subset S\}.$$
By Lemma \ref{c1} $V^\infty_S(X)$ is a subalgebra of $V^\infty(X)$
(without unit, unless $S=X$). If $S$ is closed then
$V^\infty_S(X)$ is a closed subalgebra in $V^\infty(X)$. Also
\begin{eqnarray}\label{ce1}
V^\infty_c(X)=\cup_{S\, compact} V^\infty_S(X)= \underset{S\,
compact}{\underset{\longrightarrow }{\lim}}V^\infty_S(X).
\end{eqnarray}
Let us equip $V^\infty_c(X)= \underset{S\,
compact}{\underset{\longrightarrow }{\lim}}V^\infty_S(X)$ with the
linear topology of inductive limit when each space $V^\infty_S(X)$
is equipped with the topology induced from $V^\infty(X)$. It is
easy to see that $V^\infty_c(X)$ is a locally convex Hausdorff
linear topological space. The identical imbedding
$V^\infty_c(X)\inj V^\infty(X)$ is continuous.

For any subset $S\subset X$ let us denote
\begin{eqnarray*}
W_{i,S}:=V^\infty_S(X)\cap W_i(X),\\
W_{i,c}:=V^\infty_c(X)\cap W_i(X).
\end{eqnarray*}
If $S$ is closed then $W_{i,S}\subset W_i(X)$ is a closed
subspace. We will need the following lemma.
\begin{lemma}\label{c2}
Let $S$ be a closed subset of $X$. Then for any $j=0,\dots,n$
$$H^i_S(X,\cw_j)=0 \mbox{ for } i>0$$
where $H^i_S$ denotes the $i$-th cohomology with support in $S$.
\end{lemma}
{\bf Proof.} The sheaf $\cw_j$ has the descending filtration
$$\cw_j\supset \cw_{j+1}\supset\dots\supset\cw_n.$$ It is enough
to show that for any $p$ $H^i_S(X,\cw_p/\cw_{p+1})=0$ for $i>0$.
Let us denote by $\co_X$ the sheaf of $C^\infty$-smooth functions
on $X$. Then $\cw_p/\cw_{p+1}$ is a sheaf of $\co_X$-modules. It
is well known (see e.g. \cite{godement}, Section 3.7 of Ch. II)
that on any smooth manifold $X$, for any sheaf $\cf$ of
$\co_X$-modules, and any closed subset $S\subset X$ one has
$$H^i_S(X,\cf)=0 \mbox{ for } i>0.$$
This implies the lemma. \qed

\begin{lemma}\label{c3}
(1) For any closed subset $S\subset X$ the canonical isomorphism
\begin{eqnarray}\label{ce2}
W_i/W_{i+1}\tilde\to C^\infty(X,Val_i^{sm}(TX))
\end{eqnarray}
induces isomorphism
$$W_{i,S}/W_{i+1,S}\tilde\to C^\infty_S(X,Val_i^{sm}(TX))$$
where $C^\infty_S$ stays for the space of infinitely smooth
sections with support in $S$.

(2) Similarly the isomorphism (\ref{ce2}) indices isomorphism
$$W_{i,c}/W_{i+1,c}\tilde \to C^\infty_c(X,Val_i^{sm}(TX)).$$
\end{lemma}
{\bf Proof.} Part (2) follows from part (1) by passing to direct
limit. Thus let us prove part (1). Equality (\ref{comp3}) implies
that we have a short exact sequence of sheaves on $X$
$$0\to \cw_{i+1}\to \cw_i\to \underline{Val}_i(TX)\to 0.$$
Taking sections with the support in $S$ we obtain the following
exact sequence
$$0\to W_{i+1,S}\to W_{i,S}\to C^\infty_S(X,Val_i^{sm}(TX))\to H^1_S(X,\cw_{i+1}).$$
But by Lemma \ref{c2} $H^1_S(X,\cw_{i+1})=0$. Hence the result
follows. \qed



\subsection{The integration functional.}\label{integration-functional}
In this subsection we are going to introduce a canonical linear
functional
$$\int \colon V^\infty_c(X)\to \CC$$
which we call the {\itshape integration functional}. With slight
oversimplification
$$\int \phi=\phi(X)$$
for any $\phi\in V^\infty_c(X)$. This definition is formally
correct if $X$ is compact. Otherwise $X\not\in \cp(X)$, and the
above definition requires an explanation.

Let us construct the integration functional formally for general
manifold $X$. First fix a compact subset $S\subset X$. Let us
choose a compact subset $S'$ with smooth boundary and such that
$S$ in contained in the interior of $S'$. Then $S'\in \cp(X)$. For
any $\phi\in V^\infty_S(X)$ define
\begin{eqnarray}\label{ce3}
\si \phi :=\phi(S').
\end{eqnarray}
\begin{lemma}\label{c4}
(1) $\si\colon V^\infty_S(X)\to \CC$ is a continuous linear
functional.

(2) For fixed $S$, the right hand side in (\ref{ce3}) is
independent of $S'$ containing $S$.
\end{lemma}
{\bf Proof.} Part (1) is obvious. Let us prove part (2). Let $S''$
be another compact subset with smooth boundary containing $S$ in
the interior. We have to show that $\phi(S')=\phi(S'')$. Choosing
a larger subset if necessary one may assume that $S'$ is contained
in the interior of $S''$. Then
$$\phi(S'')=\phi(S')+\phi(\overline{S''\backslash S'})-\phi(\pt S')=\phi(S')$$
where the last equality is due to the fact that $\supp\,
\phi\subset S\subset int S'$. \qed

As in the proof of Lemma \ref{c4} it is easy to see that if
$S_1\subset S_2$ then the restriction of $\sit$ to
$V^\infty_{S_1}(X)$ is equal to $\sio$. Thus we obtain a {\itshape
continuous} linear functional
$$\int\colon V^\infty_c(X)\to \CC.$$
\begin{remark}
The space of smooth compactly supported densities is a subspace of
$V^\infty_c(X)$; it is equal to $W_{n,c}$. The restriction of the
above constructed integration functional $\int$ to this subspace
coincides with the usual integration of densities.
\end{remark}

\section{The selfduality property of
valuations.}\label{selfduality-section} The goal of this section
is to establish the selfduality property of valuations (Theorem
\ref{s1}, Subsection \ref{selfduality-theorem}). Subsection
\ref{partition-unity} contains a technical result on partition of
unity in valuations.
\subsection{The selfduality property.}\label{selfduality-theorem} Probably the most interesting property of
the multiplicative structure on valuations is Theorem \ref{s1}
below. Its prove heavily uses the Irreducibility Theorem for
translation invariant valuations from \cite{alesker-gafa}.
\begin{theorem}\label{s1}
Consider the bilinear form
$$V^\infty(X)\times V^\infty_c(X)\to \CC$$
given by $(\phi,\psi)\mapsto \int \phi\cdot \psi$.

This bilinear form is a perfect pairing. More precisely the
induced map
$$V^\infty(X)\to (V^\infty_c(X))^*$$
is injective and has a dense image with respect to the weak
topology on $(V^\infty_c(X))^*$.
\end{theorem}
Theorem \ref{s1} follows from the next more precise statement by
application of the Hahn-Banach theorem.
\begin{theorem}\label{s1.1}
(1) For any $\phi\in W_i\backslash W_{i+1}$ there exists $\psi\in
W_{n-i,c}$ such that $\int \phi\cdot \psi \ne 0$.

(2) Similarly for any $\phi\in W_{i,c}\backslash W_{i+1,c}$ there
exists $\psi\in W_{n-i}$ such that $\int \phi\cdot \psi \ne 0$.
\end{theorem}

{\bf Proof.} The proves of these two statements are very similar.
Thus let us prove only the first one. Let $\phi\in W_i\backslash
W_{i+1}$. Let us denote by $\tilde\phi$ the image of $\phi$ is
$W_i/W_{i+1}=C^\infty(X,Val_i^{sm}(TX))$. Thus $\tilde \phi\ne 0$.
We will show that there exists $\psi\in W_{n-i,c}$ such that
$\int\phi\cdot \psi\ne 0$. Since $W_{i+1}\cdot W_{n-i}=0$ and
$W_i\cdot W_{n-i+1}=0$, the product $\phi\cdot \psi$ depends only
on $\tilde\phi$ and on the image $\tilde\psi$ of $\psi$ in
$W_{n-i,c}/W_{n-i+1,c}= C^\infty_c(X,Val_{n-i}^{sm}(TX))$ (where
the last equality is due to Lemma \ref{c3}(2).

Thus it is enough to show that for any $\tilde\phi\in
C^\infty(X,Val_i^{sm}(TX))$ there exists $\tilde\psi\in
C^\infty_c(X,Val_{n-i}^{sm}(TX))$ such that
$$\int_X \tilde \phi\cdot\tilde\psi \ne 0$$
where the product $\tilde\phi\cdot \tilde\psi$ is understood
pointwise in the tangent space of each point, $\tilde\phi\cdot
\tilde\psi\in
C^\infty_c(X,Val_n^{sm}(TX))=C^\infty_c(X,|\ome_X|)$, and the
integration is understood in the sense of the usual integration of
densities.

Let us fix a point $x_0\in X$ such that $\tilde\phi(x_0)\ne 0$. By
the Poincar\'e duality for the translation invariant valuations
(Theorem 0.8 in \cite{alesker-poly}) there exists $\xi_0\in
Val_{n-i}^{sm}(T_{x_0}X)$ such that $\tilde\phi(x_0)\cdot \xi_0\ne
0$. Let $\xi\in C^\infty(X,Val_{n-i}^{sm}(TX))$ be a section such
that $\xi(x_0)=\xi_0$.

Consider the following $C^\infty$-smooth density on $X$
$$\tau:=\tilde\phi\cdot \xi.$$
Thus $\tau(x_0)\ne 0$. Hence we can find a smooth compactly
supported function $\delta\in C^\infty_c(X)$ such that
$\int_X\tau\cdot \delta \ne 0$. Take $\tilde \psi:=\xi\cdot
\delta$. Then
$$\int\tilde\phi\cdot \tilde\psi=\int_X\tau\cdot \delta\ne 0.$$
\qed

From Theorem \ref{s1.1} we immediately deduce the following
corollary.
\begin{corollary}\label{s1.2}
\begin{eqnarray*}
W_i=\{\phi\in V^\infty(X)|\, \int \phi\cdot \psi =0 \mbox{ for any
} \psi\in W_{n-i+1,c}\},\\
W_{i,c}=\{\phi\in V^\infty_c(X)|\, \int \phi\cdot \psi =0 \mbox{
for any } \psi\in W_{n-i+1}\}.
\end{eqnarray*}
\end{corollary}

\subsection{Partition of unity in
valuations.}\label{partition-unity}
\begin{proposition}\label{partition}
Let $\{U_\alp\}_\alp$ be a locally finite open covering of a
manifold $X$. Then there exist $\{\phi_\alp\}_\alp\subset
V^\infty(X)$ such that
$$\supp\, (\phi_\alp)\subset U_\alp  \mbox{ and }
\sum_\alp\phi_\alp \equiv\chi$$ where the sum is locally finite,
and $\chi$ denotes the Euler characteristic.
\end{proposition}
Proposition \ref{partition} is an immediate consequence of the
fact that the sheaf $\cv^\infty_X$ of smooth valuations is soft
(by Proposition 3.1.8 of \cite{part2}) and the following general
result.
\begin{proposition}[\cite{godement}, Theorem 3.6.1, Ch.
II]\label{partun} Let $X$ be a paracompact topological space. Let
$\{U_i\}_{i\in I}$ be a locally finite open covering of $X$. Let
$\cl$ be a soft sheaf over $X$. Then for any section $s\in \cl(X)$
there exists a collection of sections $\{s_i\}_{i\in I}\subset
\cl(X)$ such that

(1) $\supp s_i\subset U_i$;

(2) the family of subsets $\{\supp s_i\}_{i\in I}$ is locally
finite;

(3) $s=\sum_{i\in I} s_i$.
\end{proposition}

\section{Generalized
valuations.}\label{generalized-valuations-section} In this section
we introduce and study the space $\vmi$ of {\itshape generalized}
valuations. It is defined in Subsection \ref{gen-val}. In
Subsection \ref{sheaf-prop-gener} it is shown that generalized
valuations form naturally a sheaf on $X$; it is a sheaf of modules
over the sheaf of algebras of smooth valuations. In Subsection
\ref{filtration} a canonical filtration on generalized valuations
is introduced and studied; it extends in a sense the canonical
filtration on smooth valuations. In Subsection
\ref{euler-verdier-generalized} we extend the Euler-Verdier
involution from smooth valuations to generalized ones.
\subsection{The space of generalized valuations.}\label{gen-val}
\begin{definition}\label{2-1}
Define the space of {\itshape generalized valuations} by
$$V^{-\infty}(X):= (V_c^\infty(X))^*$$
equipped with the usual weak topology on the dual space.
\end{definition}
\begin{remark}
It is important to observe that by Theorem \ref{s1} we have a
canonical imbedding $$V^\infty(X)\inj V^{-\infty}(X)$$ with the
image dense in the weak topology. Thus we can consider the space
of generalized valuations as a completion of the space of smooth
compactly supported valuations with respect to the weak topology.
\end{remark}

Let us describe on $V^{-\infty}(X)$ the canonical structure of
$V^\infty(X)$-module. Let $\xi\in V^\infty(X),\, \psi\in
V^{-\infty}(X)$. Define their product $\xi\cdot \psi$ by
$$<\xi\cdot \psi,\phi>=<\phi,\xi\cdot \phi>$$
for any $\phi\in V^\infty_c(X)$. Clearly this defines a map
$$\mu\colon V^\infty(X)\times V^{-\infty}(X)\to V^{-\infty}(X).$$
\begin{proposition}\label{2-2}
The map $\mu$ is a separately continuous bilinear map. It defines
a structure of $V^\infty(X)$-module on $V^{-\infty}(X)$. Moreover
$V^\infty(X)$ is a submodule of $V^{-\infty}(X)$, and the induced
structure of $V^\infty(X)$-module on it is the standard one.
\end{proposition}
{\bf Proof.} The bilinearity is obvious from the definition. Let
us check the continuity. We have to check that for any $\phi\in
V_c^\infty(X)$ the map
$$V^\infty(X)\times V^{-\infty}(X)\to \CC$$
given by $(\xi,\psi)\mapsto <\psi,\xi\cdot \phi>$ is separately
continuous. But this is an immediate consequence of the continuity
of the map $V^\infty(X)\to V^\infty_c(X)$ given by $\xi\mapsto
\xi\cdot \phi$ and separate continuity of the canonical pairing
$V^\infty_c(X)\times V^{-\infty}(X)\to \CC$.

Let us check now that the above map $\mu\colon V^\infty(X)\times
V^{-\infty}(X)\to V^{-\infty}(X)$ defines the standard
$V^\infty(X)$-module structure on $V^\infty(X)\inj
V^{-\infty}(X)$. Namely we have to show that for $\xi,\psi\in
V^\infty(X)$ one has $\mu(\xi,\psi)=\xi\cdot \psi$ where the last
product is understood in the usual sense. Let $\phi\in
V_c^{\infty}(X)$. Then we have
\begin{eqnarray*}
<\mu(\xi,\psi),\phi>=<\psi,\xi\cdot \phi>=\\
\int\psi\cdot (\xi\cdot \phi)=\int(\xi\cdot \psi)\cdot
\phi=<\xi\cdot\psi,\phi>.
\end{eqnarray*}
Hence $\mu(\xi,\phi)=\xi\cdot \psi$.

Since $V^\infty(X)$ is dense in $V^{-\infty}(X)$ and $\mu$ is
continuous it follows that $\mu$ defines $V^\infty(X)$-module
structure on $V^{-\infty}(X)$. \qed

\subsection{The sheaf property of generalized valuations.}\label{sheaf-prop-gener} In this
subsection we describe the canonical sheaf structure on
generalized valuations.

First observe that for two open subsets $U_1\subset U_2$ of a
manifold $X$ we have the identity imbedding
\begin{eqnarray}\label{2-4.5}
V^\infty_c(U_1)\inj V^\infty_c(U_2).
\end{eqnarray}
 Hence by duality we have a {\itshape continuous}
map
\begin{eqnarray}\label{2-5}
V^{-\infty}(U_2)\to V^{-\infty}(U_1).
\end{eqnarray}

\begin{lemma}\label{2-6}
The map (\ref{2-5}) being restricted to $V^\infty(U_2)\subset
V^{-\infty}(U_2)$ coincides with the usual restriction map
$V^\infty(U_2)\to V^\infty(U_1)$.
\end{lemma}
{\bf Proof.} Let us denote temporarily the imbedding (\ref{2-4.5})
by $\tau$, and its dual (\ref{2-5}) by $\tau^*$.  Let $\phi\in
V^\infty(U_2)$. Then for any $\psi\in V^\infty_c(U_1)$ one has
\begin{eqnarray*}
<\tau^*(\phi),\psi>=(\phi\cdot \tau(\psi))(U_2)=(\phi|_{U_1}\cdot
\psi)(U_1)=<\phi|_{U_1},\psi>.
\end{eqnarray*}
Hence $\tau^*(\phi)=\phi|_{U_1}$. \qed

\begin{proposition}\label{2-7}
The assignment
$$U\mapsto V^{-\infty}(U)$$
to any open subset $U\subset X$ with the above restriction maps
defines a sheaf on $X$ denoted by $\cv^{-\infty}_X$.
\end{proposition}
\begin{remark}\label{2-8} Given this proposition, it is clear that
$\cv^\infty_X$ is a subsheaf of $\cv^{-\infty}_X$.
\end{remark}
{\bf Proof} of Proposition \ref{2-7}. Let $\{U_\alp\}_\alp$ be an
open covering of an open subset $U$. Let $\phi\in V^{-\infty}(V)$
such that $\phi|_{U_\alp}= 0$ for any $\alp$. Let us show that
$\phi =0$. Replacing $\{U_\alp\}$ by a refinement we may assume
that $\{U_\alp\}$ is locally finite. Let us choose a partition of
unity $\{\phi_\alp\}$ subordinate to this covering using
Proposition \ref{partition}. For any $\psi\in V_c^\infty(X)$ we
have
\begin{eqnarray*}
<\phi,\psi>=<\phi,\sum_\alp \phi_\alp\cdot \psi>=\\
\sum_\alp<\phi,\phi_\alp\cdot
\psi>=\sum_\alp<\phi|_{U_\alp},(\phi_\alp\cdot \psi)_{U_\alp}>=0.
\end{eqnarray*}
Hence $\phi=0$.

Now let us assume that we are given an open covering
$\{U_\alp\}_\alp$ of an open subset $U\subset X$, and  for any
$\alp$ we are given a generalized valuation $\psi_\alp\in
V^{-\infty}(U_\alp)$ such that $\psi_\alp|_{U_\alp\cap
U_\beta}=\psi_\beta|_{U_\alp\cap U_\beta}$ for any $\alp, \beta$.
Let us show that there exists $\psi\in V^{-\infty}(U)$ such that
$\psi|_{U_\alp}=\psi_\alp$. Again by choosing a refinement we may
assume that the covering $\{U_\alp\}$ is locally finite. Let us
fix a partition of unity $\{\phi_\alp\}$ subordinate to it. Define
$\psi$ by
\begin{eqnarray*}
<\psi,\phi>:=\sum_\alp <\psi_\alp,(\phi_\alp\cdot \phi)|_{U_\alp}>
\end{eqnarray*}
for any $\phi\in V^\infty_c(U)$. It is easy to see that $\psi\in
V^{-\infty}(U)$ and $\psi|_{U_\alp}=\psi_\alp$. \qed

\begin{proposition}\label{2-9}
Being equipped with the above restriction maps and the defined
above product of generalized valuations by smooth ones,
$\cv^{-\infty}_X$ is a sheaf of $\cv^\infty_X$-modules.
\end{proposition}
{\bf Proof.} For an open subset $U\subset X$ let us denote by
$$\mu_U\colon V^\infty(U)\times V^{-\infty}(U)\to
V^{-\infty}(U)$$ the canonical product. We have to check that for
any open subsets $U\subset V\subset X$, any $\xi\in V^\infty(V),\,
\psi\in V^{-\infty}(V)$ one has
\begin{eqnarray}\label{2-10}
\left(\mu_V(\xi,\psi)\right)|_U=\mu_U(\xi|_U,\psi|_U).
\end{eqnarray}
Let $\phi\in V^\infty_c(U)$. Let us denote the identity imbedding
$V^\infty_c(U)\inj V^\infty_c(V)$ by $\tau$.

Then we have
\begin{eqnarray*}
<(\mu_V(\xi,\psi))|_U,\phi>=<\mu_V(\xi,\psi),\tau(\phi)>=
<\psi,\xi\cdot \tau(\phi)>=\\
<\psi,\tau(\xi|_U\cdot \phi)>=<\psi|_U,\xi|_U\cdot
\phi>=<\mu_U(\xi|_U,\psi|_U),\phi>.
\end{eqnarray*}
Hence (\ref{2-10}) follows. \qed

\subsection{Filtration on generalized
valuations.}\label{filtration}
\begin{definition}\label{f1}
Define $W_i(V^{-\infty}(X))$ to be the closure of $W_i(X)(\subset
V^\infty(X)\subset V^{-\infty}(X))$ in the space $\vmi$  with
respect to the weak topology.
\end{definition}
Clearly one has
\begin{eqnarray*}
V^{-\infty}(X)=W_0(\vmi)\supset W_1(\vmi)\supset \dots \supset
W_n(\vmi).
\end{eqnarray*}
In this subsection we will also use the following notation. The
subspace  $W_i(X)$ of $V^\infty(X)$ will also be denoted by
$W_i(\vi)$. Set
\begin{eqnarray*}
W_i(\vic):=W_i(\vi)\cap \vic,\\
W_i(\vmic):=W_i(\vmi)\cap \vmic.
\end{eqnarray*}
It is easy to see (using the separate continuity of the product
$V^\infty(X)\times V^{-\infty}(X)\to V^{-\infty}(X)$) that
\begin{eqnarray}\label{inclus}
W_i(V^\infty(X))\cdot W_j(V^{-\infty}(X))\subset
W_{i+j}(V^{-\infty}(X)).
\end{eqnarray}
\begin{proposition}\label{f2}
For any $i=0,1,\dots,n$
\begin{eqnarray}\label{f3}
W_i(\vmi)=\{\phi\in \vmi|\, <\phi,\psi>=0 \mbox{ for any } \psi\in
W_{n-i+1}(\vic)\}.
\end{eqnarray}
\end{proposition}
{\bf Proof.} Let us denote by $W'_i(\vmi)$ the space in the right
hand side of (\ref{f3}). The equality (\ref{inclus}) implies that
$$W_i(\vmi)\subset W'_i(\vmi).$$
Let us prove the converse inclusion. Let us assume in the contrary
that there exists $\psi\in W'_i(\vmi)\backslash W_i(\vmi)$. Since
$W_i(\vmi)$ is a closed subspace of $\vmi$ in the weak topology,
the Hahn-Banach theorem implies that there exists $\phi\in \vic$
such that $<\psi,\phi>\ne 0$ and for any $\xi\in W_i(\vmi)$
$$<\xi,\phi>= 0.$$
Since $W_i(\vi)\subset W_i(\vmi)$ Corollary \ref{s1.2} implies
that $\phi\in W_{n-i+1}(\vic)$. But then (\ref{inclus}) implies
that $<\psi,\phi>=0$. This is a contradiction. \qed
\begin{corollary}\label{f4}
$$W_i(\vmi)\cap \vi= W_i(\vi).$$
\end{corollary}
{\bf Proof.} This immediately follows from Proposition \ref{f2}
and Corollary \ref{s1.2}. \qed

For a subset $S\subset X$ let us denote by $V^{-\infty}_S(X)$ the
space of generalized valuations with support contained in $S$.
Clearly $V^{-\infty}_S(X)$ is a $\vi$-submodule of $\vmi$. If $S$
is a closed subset of $X$ then $V^{-\infty}_S(X)$ is a closed
subspace of $\vmi$ in the weak topology. It is easy to see that
$$\vmic =\underset{S\, compact}{\underset{\longrightarrow }{\lim}} V^{-\infty}_S(X).$$

Let us equip $\vmic$ with the topology of inductive limit when
each of $V^{-\infty}_S(X)$ is equipped with the topology induced
from $\vmi$. Then $\vmic$ is a locally convex Hausdorff linear
topological space.
\begin{proposition}\label{f5}
For any $i=0,1,\dots,n,$ the space $W_i(\vic)$ is dense in
$W_i(\vmic)$ in the above topology of inductive limit.
\end{proposition}
{\bf Proof.} Fix $\phi\in W_i(\vmic)$. Set $S:=\supp \phi$ be the
support of $\phi$. $S$ is a compact set. Let $U$ be an open
relatively compact neighborhood of $S$. Since the sheaf
$\cv^\infty_X$ of smooth valuations is soft (by Proposition 3.1.8
of \cite{part2} ), there exists $\alp\in V^\infty(X)$ such that
$\alp$ is equal to the Euler characteristic $\chi$ in a
neighborhood of $S$, and $\alp|_{X\backslash U}\equiv 0$.

Since $W_i(\vi)$ is dense in $W_i(\vmi)$ in the weak topology,
there exists a net $\{\phi_\lam\}\subset W_i(\vi)$ which converges
to $\phi$ in the weak topology. But then $\{\alp\cdot
\phi_\lam\}\subset W_i(V^\infty_U(X))$, and $\{\alp\cdot
\phi_\lam\}$ converges to $\alp\cdot\phi=\phi$ in $\vmic$. \qed

\begin{proposition}\label{f6}
Let $n=\dim X$ as previously. Then there exists a canonical
isomorphism of linear topological spaces
$$W_n(\vmi)=C^{-\infty}(X,|\ome_X|).$$
\end{proposition}
{\bf Proof.} By Proposition \ref{f2} one has
\begin{eqnarray*}
W_n(\vmi)=W_1(\vic)^\perp=(\vic/ W_1(\vic))^*=\\
(C_c^\infty(X))^*=C^{-\infty}(X,|\ome_X|)
\end{eqnarray*}
where the third equality is due to Lemma \ref{c3}(2).\qed

\begin{proposition}\label{f7}
There exists a canonical isomorphism of linear topological spaces
$$\vmi/W_1(\vmi)=C^{-\infty}(X).$$
\end{proposition}
{\bf Proof.} Using Proposition \ref{f2} one has
\begin{eqnarray*}
\vmi/W_1(\vmi)=\vic^*/W_n(\vic)^\perp=\\
W_n(\vic)^*=(C^\infty_c(X,|\ome_X|))^*= C^{-\infty}(X).
\end{eqnarray*}
\qed

Remind that by Proposition \ref{f5} $\vic$ is dense in $\vmic$ (in
the topology of inductive limit).
\begin{proposition}\label{f8}
The integration functional
$$\int\colon \vic\to \CC$$
extends uniquely by continuity to the functional
$$\int\colon \vmic\to \CC.$$
\end{proposition}
{\bf Proof.} First observe that for any $\alp\in \vic$ the
functional $\vi\to \CC$ given by $\phi\mapsto \int \alp\cdot \phi$
extends (uniquely) by continuity in the weak topology to $\vmi$.
Indeed this extension is given by $\psi\mapsto <\psi,\alp>$. Let
us denote this functional by $\hat\alp$. Thus $\hat\alp\colon
\vmi\to \CC$ is a continuous functional.

Let us fix an arbitrary compact subset $S\subset X$. Let us fix a
smooth compactly supported valuation $\alp\in \vic$ such that
$\alp$ equals to the Euler characteristic in a neighborhood of
$S$. Consider the corresponding continuous linear functional
$\hat\alp\colon \vmi\to \CC$. We claim that the restriction of
$\hat \alp$ to $V^{-\infty}_S(X)$ is the desired extension of the
integration functional to $V^{-\infty}_S(X)$.

To check it let us fix a compact neighborhood $S'$ of $S$ such
that the restriction of $\alp$ to $S'$ is still equal to the Euler
characteristic. By (the proof of) Proposition \ref{f5} every
valuation from $V^{-\infty}_S(X)$ can be approximated in the weak
topology by a net from $V^\infty_{S'}(X)$. Hence it is enough to
check that for any $\phi\in V^\infty_{S'}(X)$ one has
$$\int \alp\cdot \phi=\int \phi.$$ But this is obvious since
$$\phi\cdot (\alp-\chi) \equiv 0.$$ \qed





\begin{lemma}\label{f8.5}
Let $\{\zeta_\lam\}_{\lam\in \Lam}\subset \vi$ be a net such that
for any compact subset $K\subset X$ there exists $\lam_K\in \Lam$
such that for all $\lam\geq \lam_K$
$$(\supp \zeta_K)\cap K=\emptyset.$$
Then $$\lim_\Lam \zeta_\lam=0 \mbox{ in } \vi.$$
\end{lemma}
{\bf Proof.} Consider the map
$$T_i\colon \tilde C^\infty(T^*X,W_i(\Ome^n\otimes p^*o))\to
W_i(\vi)$$ given by the integration with respect to the
characteristic cycle. By Proposition 3.1.3 of \cite{part2} $T_i$
is an epimorphism. By the definition of the topology on $\vi$ (see
Subsection 3.2 of \cite{part2}) $T_i$ is a continuous map. Hence
it is enough to show that for any compact subset $K\subset X$
there exists $\lam_K\in \Lam$ such that for any $\lam\geq\lam_K$
there exists  $\eta_\lam\in \tilde C^\infty(T^*X,\Ome^n\otimes
p^*o)$ satisfying

(i) $T_i(\eta_\lam)=\zeta_\lam$;

(ii) $\eta_\lam$ vanishes in a neighborhood of $p^{-1}(K)$.

Indeed then we would have
$$\lim_{\Lam}\zeta_\lam=\lim_{\Lam}T_i(\eta_\lam)=0.$$

For the rest of the proof of the lemma let fix a compact subset
$K\subset X$. As in Section \ref{technical-lemma} consider the
sheaves on $X$
\begin{eqnarray*}
\cw_i(U)=W_i(V^\infty(U)),\\ \cw_i'(U)=\tilde
C^\infty(T^*U,W_i(\Ome^n\otimes p^*o))
\end{eqnarray*}
for any open subset $U\subset X$. The integration over the the
characteristic cycle gives a morphism of sheaves
\begin{eqnarray*}
T_i\colon \cw'_i\to \cw_i
\end{eqnarray*}
which is an epimorphism (we denote this morphism by the same
symbol $T_i$). Set again $\ck_i:=Ker T_i$. It was shown in
\cite{part2} (see the proof of Proposition 3.1.9) that the sheaves
$\ck_j/\ck_{j+1}$ are isomorphic to the sheaves of $\co_X$-modules
where $\co_X$ denotes the sheaf of $C^\infty$-smooth functions on
$X$. By Section 3.7 of Ch. II of \cite{godement} the sheaves
$\ck_j/\ck_{j+1}$ are soft for any $j$. Hence for any closed
subset $Z\subset X$ the positive cohomology groups with support in
$Z$ vanish:
$$H^i_Z(X,\ck_j/\ck_{j+1})=0 \mbox{ for } i>0.$$ Using the long
exact sequence we get
$$H^i_Z(X,\ck_j)=0 \mbox{ for } i>0 \mbox{ and any } j.$$
Consider the short exact sequence of sheaves
$$0\to\ck_i\to \cw_i'\to\cw_i\to 0.$$
From the long exact sequence we obtain
\begin{eqnarray}\label{zz1}
H^0_Z(X,\cw_i')\to H^0_Z(X,\cw_i)\to H^1_Z(X,\ck_i)=0
\end{eqnarray}
for any closed subset $Z\subset X$. Namely the map
\begin{eqnarray}\label{zz2}
H^0_Z(X,\cw_i')\to H^0_Z(X,\cw_i)
\end{eqnarray}
is surjective.

Let us choose $Z$ as follows. Let $U$ be an open relatively
compact neighborhood of $K$. Set $Z:=X\backslash U$. There exists
$\lam_0\in \Lam$ such that for any $\lam\geq \lam_0$ one has
$(\supp \zeta_\lam)\cap U=\emptyset$. Then clearly $\zeta_\lam\in
H^0_Z(X,\cw_i)$ for $\lam\geq \lam_0$. The surjectivity of the map
(\ref{zz2}) implies the lemma. \qed



\begin{lemma}\label{f9}
(1) For any $i=0,1,\dots,n$ the space $W_i(V^\infty_c(X))$ is
dense in $W_i(\vi)$.

(2) For any $i=0,1,\dots,n$ the space $W_i(\vmic)$ is dense in
$W_i(\vmi)$.
\end{lemma}
{\bf Proof.}  Let us prove first part (1). For any compact subset
$K\subset X$ let us choose a compactly supported valuation
$\tau_K\in \vic$ such $\tau_K$ is equal to the Euler
characteristic $\chi$ in a neighborhood of $K$. Let $\psi\in
W_i(\vi)$. It is enough to show that
$$\lim_{K \,compact}(\tau_K\cdot \psi)=\psi
\mbox{ in } \vi.$$ Let us denote $\zeta_K:=(\tau_K-\chi)\cdot
\psi$. Clearly $\zeta_K$ vanishes in a neighborhood of $K$. By
Lemma \ref{f8.5} $\lim_{K \, compact}\zeta_K=0$.



Let us prove part (2). Fix $\psi\in W_i(\vmi)$. For any compact
subset $K\subset X$ let us fix  a compactly supported smooth
valuation $\tau_K\in \vic$ which is equal to the Euler
characteristic $\chi$ in a neighborhood of $K$. Let
$\psi_K:=\tau_K\cdot \psi$. Clearly $\psi_K\in W_i(\vmic)$. It
suffices to show that
$$\lim_{K\, compact} \psi_K=\psi \mbox{ in
} \vmi.$$ Let $\phi\in \vic$. We have to show that
$$\lim_{K\, compact}<\psi_K,\phi>=<\psi,\phi>.$$ We have
\begin{eqnarray*}
\lim_{K\, compact}<\psi_K,\phi>=\lim_{K\,
compact}<\psi,\tau_K\cdot \phi>=\\<\psi,\lim_{K\,
compact}(\tau_K\cdot \phi)>=<\psi,\phi>.
\end{eqnarray*}
Part (2) is proved too. \qed


Let us observe now that the bilinear map $\vmic\times \vi\to \CC$
given by
$$(\psi,\phi)\mapsto \int\phi\cdot \psi$$
is separately continuous. Hence it defines a continuous map
$$\theta\colon \vmic \to \vi^*$$
where $\vi^*$ is equipped with the weak topology, and a continuous
map
$$\theta'\colon \vi\to \vmic^*$$
where $\vmic^*$ is equipped with the weak topology.
\begin{proposition}\label{f10}
The maps $\theta$ and $\theta'$ are isomorphisms of linear spaces.
\end{proposition}
{\bf Proof.} First observe that if the manifold $X$ is compact
then the result follows immediately from the definitions. Let us
assume that $X$ is not compact.

First let us check that $\theta$ is injective. Assume that
$\psi\in Ker \theta$. Since $\psi\in \vmic\subset \vmi=\vic^*$,
then for any $\phi\in \vic$ one has $<\psi,\phi>=0$. Hence
$\psi=0$.

Let us check now that $\theta$ is onto. Let $\zeta\in \vi^*$.
Since the identity imbedding $\vic\hookrightarrow\vi$ is
continuous, the restriction $\tilde\zeta$ of $\zeta$ to $\vic$ is
a continuous functional on $\vic$, i.e. belongs to $\vic^*=\vmi$.
Let us show that $\supp \tilde \zeta$ is compact, i.e.
$\tilde\zeta\in \vmic$. Assume in the contrary that $\supp \tilde
\zeta$ is not compact. It means that for any compact subset
$K\subset X$ there exists a valuation $\phi\in \vic$ with $\supp
\phi \cap K=\emptyset$ such that $<\tilde\zeta,\phi>\ne 0$. Since
we have assumed that $X$ is not compact we can construct an open
covering $\{U_\alp\}_\alp$ of $X$ which does not have a finite
subcovering. Since any manifold is paracompact (by definition) and
locally compact, by choosing a refinement if necessary we may
assume that this covering is locally finite and any $U_\alp$ is
relatively compact. Let us choose $U_{\alp_1}$ so that $\supp
\tilde\zeta\cap U_{\alp_1}\ne \emptyset$. Denote $K_1:=\bar
U_{\alp_1}$. Assume we have constructed compact sets
$K_1,\dots,K_{N-1}$ with the following properties:

1) for each $i=1,\dots,N-1$ there exists $\alp_i$ such that
$K_i=\bar U_{\alp_i}$;

2) the interior of $K_i$ intersects $\supp \tilde\zeta$
non-trivially for each $i=1,\dots,N-1$;

3) $K_i\cap K_j=\emptyset$ for $1\leq i\ne j\leq N-1$.

Let us construct $K_N$ such that the sequence of sets
$K_1,\dots,K_{N-1},K_N$ has the same properties. Let us fix an
open relatively compact neighborhood $T$ of the set
$\cup_{i=1}^{N-1}K_i$. Since the covering $\{U_\alp\}$ is locally
finite, and $\supp \tilde\zeta$ is not compact, there exists
$\alp_N$ such that $U_{\alp_N}\cap T=\emptyset$ and
$U_{\alp_N}\cap \supp \tilde\zeta \ne\emptyset$. Set $K_N:=\bar
U_{\alp_N}$. Then $K_N\cap (\cup_{i=1}^{N-1}K_i)=\emptyset$ and
$K_N\cap \supp\tilde\zeta\ne \emptyset$. By induction we obtain an
infinite sequence of pairwise disjoint compact sets $\{K_N\}_{N\in
\NN}$ with non-empty interiors such that $int K_N\cap \supp
\tilde\zeta\ne \emptyset$ for any $N\in \NN$.

Since $int K_N\cap\supp \tilde\zeta\ne \emptyset$ we can choose a
valuation $\phi_N\in V^\infty(X)$ with $\supp \phi_N\subset int
K_N$ and such that $<\tilde\zeta,\phi_N>=1$. Let us define
$$\phi:=\sum_{N=1}^\infty\phi_N.$$
This series converges in $\vi$ by Lemma \ref{f8.5}. Then
$$<\zeta,\phi>=\lim_{N\to \infty}<\zeta,\sum_{n=1}^N\phi_n>=\lim_{N\to
\infty}N=\infty.$$ This is a contradiction. Hence we have shown
that $\supp \tilde \zeta$ is compact.

Let us show that $\zeta=\theta(\tilde \zeta)$. For any $\phi\in
\vic$ we have
$$<\zeta,\phi>=<\tilde\zeta,\phi>.$$ Hence
$<\zeta,\phi>=<\theta(\tilde\zeta),\phi>$ for any $\phi\in \vic$.
But by Lemma \ref{f9} $\vic$ is dense in $\vi$. Hence
$\zeta=\theta(\tilde\zeta)$. Thus we have shown that $\theta\colon
\vmic\to \vi^*$ is an isomorphism of linear spaces.

Let us show that $\theta'$ is an isomorphism of linear spaces.
First let us check that $\theta'$ is injective. Assume that
$\psi\in Ker \theta'$. Since $\vic\subset \vmic$ then for any
$\phi\in \vic$ one has
$$\int \phi\cdot \psi =0.$$
By the Selfduality Property (Theorem \ref{s1}) $\psi\equiv 0$.

Let us show that $\theta'$ is surjective. Let $\zeta\in \vmic^*$.
For any compact subset $K\subset X$ let us fix a compactly
supported valuation $\gamma_K\in \vic$ such that the restriction
of $\gamma_K$ to a neighborhood of $K$ is equal to the Euler
characteristic $\chi$. Consider the linear functional
$$\zeta_K\colon \vmi\to \CC$$
defined by $\zeta_K(\phi)=\zeta(\gamma_K\cdot \phi)$. It is easy
to see that $\zeta_K$ is a continuous functional on $\vmi$
equipped with the weak topology. Hence $\zeta_K\in \vmic^*=\vic$.
It is also clear that if $K_1\subset K_2$ then the restriction of
$\zeta_{K_2}$ to $K_1$ is equal to $\zeta_{K_1}$. Taking limit
over all compact subsets of $X$ we get a smooth valuation on $X$
denoted by $\tilde \zeta$. Clearly the restriction of $\tilde
\zeta$ to any compact subset $K\subset X$ is equal to $\zeta_K$.
Then evidently $\zeta=\theta'(\tilde\zeta)$. \qed


\subsection{The Euler-Verdier involution on generalized
valuations.}\label{euler-verdier-generalized} We are going to
extend the Euler-Verdier involution from smooth valuations to
generalized ones.
\begin{theorem}\label{ev-gener}
(i) There exists unique continuous in the weak topology linear map
\begin{eqnarray}\label{evg1}
\sigma\colon V^{-\infty}(X)\to V^{-\infty}(X)
\end{eqnarray}
such that the restriction of it to $\vi$ is the Euler-Verdier
involution on smooth valuations.

(ii) $\sigma^2=Id.$

(iii) $\sigma$ commutes with the restrictions to open subsets of
$X$, and thus induces an involution of the sheaf $\cv^{-\infty}_X$
of generalized valuations (defined in Subsection
\ref{sheaf-prop-gener}).

(iv) $\sigma(W_i(\vmi))=W_i(\vmi)$ for any $i=0,1,\dots,n$.

(v) For any $\phi\in \vi,\, \xi\in \vmi$ one has
\begin{eqnarray}\label{evg2}
\sigma(\phi\cdot\xi)=\sigma(\phi)\cdot\sigma(\xi).
\end{eqnarray}
\end{theorem}
{\bf Proof.} Let us prove first part (i) The uniqueness is obvious
since $\vi$ is dense in $\vmi$ in the weak topology. Let us probe
the existence.

We have the Euler-Verdier involution on smooth valuations
$$\sigma\colon \vi\to \vi.$$
Since this map commutes with restrictions to open subsets of $X$,
it preserves support of a smooth valuation. Hence
$\sigma\colon\vic\to \vic$ is a continuous operator (with respect
to the topology of inductive limit on $\vic$). Consider the dual
operator
$$\sigma^*\colon \vmi\to \vmi.$$
$\sigma^*$ is continuous in the weak topology. Let us show that
the restriction of $(-1)^n\sigma^*$ to smooth valuations coincides
with the Euler-Verdier involution on $\vi$. This will finish the
proof of part (i) since $(-1)^n\sigma^*$ is the operator we need
(which will be denoted again by $\sigma$).

Let $\psi\in \vi\subset \vmi$. It is enough to show that for any
$\phi\in \vic$ one has
\begin{eqnarray*}
<(-1)^n\sigma^*\psi,\phi>=<\sigma\psi,\phi>.
\end{eqnarray*}
Using the automorphism property of the Euler-Verdier involution on
smooth valuations (Theorem \ref{euler-verdier}) we have
\begin{eqnarray*}
<(-1)^n\sigma^*\psi,\phi>=(-1)^n<\psi,\sigma\phi>=(-1)^n\int\psi\cdot\sigma\phi=\\
\int\sigma(\psi\cdot\sigma\phi)=\int\sigma\psi\cdot\phi=<\sigma\psi,\phi>.
\end{eqnarray*}
Part (i) is proved. The remaining statements of the theorem follow
from the continuity and the corresponding properties of the
Euler-Verdier involution on smooth valuations. \qed

\section{Valuations on real analytic manifolds.}\label{real-analytic}
\setcounter{subsection}{0}\setcounter{theorem}{0}\setcounter{equation}{0}
The goal of this section is to make a comparison of valuations
with a more familiar space of constructible functions on a real
analytic manifold. Let us fix a real analytic manifold $X$ of
dimension $n$.

In Subsection \ref{construction-imbedding} we construct a
canonical imbedding of the space of constructible functions
$\cf(X)$ into the space of generalized valuations $\vmi$ as a
dense subspace. In Subsection \ref{comparison-filtrations} we show
that the restriction of the canonical filtration on $\vmi$ to
$\cf(X)$ is the filtration of $\cf(X)$ by codimension of the
support. In Subsection \ref{comparison-integration} it is proved
that the restriction of the integration functional on the space of
generalized valuations with compact support to the subspace
$\cf_c(X)$ of constructible functions with compact support is the
integration with respect to the Euler characteristic. In
Subsection \ref{ev-verdier} we show that the restriction of the
Euler-Verdier involution on generalized valuation to $\cf(X)$
coincides (up to a sign) with the Verdier duality operator on the
latter.

\subsection{Imbedding of constructible functions to generalized
valuations.}\label{construction-imbedding} In this subsection we
will construct a canonical $\CC$-linear map
$$\Xi\colon \cfx\to \vmi$$ and prove that it is injective and has
a dense image in the weak topology, where $\cf(X)$ is the space of
constructible functions on $X$ defined in Subsection \ref{suban}
(equality (\ref{constr-def})).

The construction of the map $\Xi$ is based on the notion of
characteristic cycle attached to an arbitrary constructible
function $f\in \cfx$ denoted by $CC(f)$. This notion was discussed
in Subsection \ref{chc}.




Note in addition that the characteristic cycle satisfies
\begin{eqnarray}\label{im2}
CC(\alp f+\beta g)=\alp CC(f)+\beta CC(g)
\end{eqnarray}
for any $\alp,\beta\in \CC$ (see \cite{kashiwara-schapira}, \S
9.7).



Now let us describe the canonical map
\begin{eqnarray}\label{im4}
\Xi\colon \cfx\to \vmi=(\vic)^*.
\end{eqnarray}
\def\coz{C^\infty_c(T^*X,\Ome^n\otimes p^*o)}
Let us denote by $\coz$ the space of $C^\infty$-sections with
{\itshape compact support} of the bundle $\Ome^n\otimes p^*o$ over
$T^*X$. By Lemma \ref{ccc4} we have the canonical continuous
epimorphism
\begin{eqnarray}\label{im5}
\coz\twoheadrightarrow \vic
\end{eqnarray}
given by
\begin{eqnarray}\label{im6}
\ome\mapsto [P\mapsto \int_{CC(P)}\ome]
\end{eqnarray}
for any $P\in \cp(X)$. For any constructible function $f\in \cfx$
let us define $\Xi(f)$ by
\begin{eqnarray}\label{im7}
<\Xi(f),\phi>=\int_{CC(f)}\ome
\end{eqnarray}
where $\ome\in \coz$ is an arbitrary lift of $\phi$. Once we show
that $\Xi(f)$ is well defined, then automatically it is a
continuous linear functional on $\vic$.

Thus it remains to check that $\Xi$ is well defined. More
explicitly, assume that $\ome\in \coz$ satisfies
\begin{eqnarray}\label{im8}
\int_{CC(P)}\ome=0
\end{eqnarray}
for any $P\in \cp(X)$. We have to check that
\begin{eqnarray}\label{im9}
\int_{CC(f)}\ome=0
\end{eqnarray}
for any constructible function $f\in \cfx$.

Let us fix such an $\ome$. By (\ref{im2}) it is enough to assume
that $f$ is the indicator function of a subanalytic subset $Q$.

Let us observe first of all that (obviously) every point $x\in X$
has a compact subanalytic neighborhood (and also an open
subanalytic neighborhood). Hence we can choose a compact
subanalytic neighborhood $S$ of the support of $\ome$. It is
enough to check that for any subanalytic subset $Q\subset S$ one
has
$$\int_{CC(Q)}\ome=0.$$

Any point $x\in X$ has a pair of subanalytic neighborhoods
$U_x\subset V_x$ such that $U_x$ is compact, $V_x$ is open, and
there exists a real analytic diffeomorphism $g_x\colon
V_x\tilde\to \RR^n$. Hence one can find a finite covering of $S$
by compact subanalytic subsets $\{U_i\}_i$, find open subanalytic
subsets $\{V_i\}_i$ with $U_i\subset V_i$, and real analytic
diffeomorphisms $f_i\colon V_i\tilde\to \RR^n$.

By the linearity of the characteristic cycle (\ref{im2}),
intersecting $Q$ with each $U_i$ we may assume that $Q$ is
relatively compact subset of $V_{i_0}$ for some $i_0$. Thus it
remains to prove the following statement.
\begin{lemma}\label{im10}
Let $\ome\in \tilde C^\infty(T^*\RR^n,\Ome^n\otimes p^*o)$
satisfies $$\int_{CC(P)}\ome=0\mbox{ for any } P\in \cp(X).$$

Then for any bounded subanalytic subset $Q\subset \RR^n$ one has
$$\int_{CC(Q)}\ome =0.$$
\end{lemma}
{\bf Proof.} 
We will reduce the proof of the lemma to Theorem 1 of
\cite{bernig-brocker}. Let us fix an orientation on $\RR^n$. Let
$\psi$ denote the restriction of $\ome$ to the zero section
$\underline{0}$ of $T^*X$. Thus $\psi\in C^\infty(\RR^n,\Ome^n)$.
Let $$q\colon T^*\RR^n\backslash \underline{0}\to
\PP_+(T^*\RR^n)$$ be the canonical projection. Let
$\tilde\ome:=q_*\ome$ be the integration of
$\ome|_{T^*\RR^n\backslash\underline{0}}$ along the fibers of $q$.
Let $a\colon \PP_+(T^*\RR^n)\to\PP_+(T^*\RR^n)$ be the canonical
(antipodal) involution described in Subsection \ref{chc}. Set
$\eta:=a^*\tilde \ome$. It is easy to see that
\begin{eqnarray*}
\int_{CC(P)}\ome=\int_{N(P)}\eta +\int_P\psi \mbox{ for any } P\in
\cp(\RR^n);\\
\int_{CC(f)}\ome=\int_{N(f)}\eta +\int_{\RR^n}f\cdot\psi \mbox{
for any } f\in \cf_c(\RR^n).
\end{eqnarray*}
Thus by assumption we get
\begin{eqnarray}\label{doc1}
\int_{N(P)}\eta+\int_P\psi=0
\end{eqnarray}
for any $P\in \cp(\RR^n)$.

It was shown in \cite{bernig-brocker}, Theorem 1, that a pair
$(\eta,\psi)$ with $\eta\in
C^\infty(\PP_+(T^*\RR^n),\Ome^{n-1}),\, \psi\in
C^\infty(\RR^n,\Ome^n)$ satisfies the equality (\ref{doc1}) for
any compact {\itshape subanalytic} subset $P$ if and only if it
satisfies the following two conditions (where $\pi\colon
\PP_+(T^*\RR^n)\to \RR^n$ is the canonical projection):
\begin{eqnarray}\label{doc2}
\int_{\pi^{-1}(x)}\eta=0 \mbox{ for any }x\in \RR^n,\\\label{doc3}
D\eta+\pi^*\psi=0
\end{eqnarray}
where $D\colon C^\infty(\PP_+(T^*\RR^n),\Ome^{n-1})\to
C^\infty(\PP_+(T^*\RR^n),\Ome^n)$ is an explicitly written
differential operator of second order (introduced by Rumin in
\cite{rumin}).

However in the proof of the "if" part of Theorem 1 in
\cite{bernig-brocker} the authors used equality (\ref{doc1}) not
for the whole class of compact subanalytic sets, but for the
subclass of compact subanalytic submanifolds with boundary. Hence
if (\ref{doc1}) is satisfied for all $P\in \cp(\RR^n)$ then the
conditions (\ref{doc2}), (\ref{doc3}) are satisfied, and hence
(\ref{doc1}) is satisfied for an {\itshape arbitrary} compact
subanalytic subset $P\subset \RR^n$ (again by Theorem 1 of
\cite{bernig-brocker}).

In order to prove our lemma it is enough to show that (\ref{doc1})
is satisfied for any bounded subanalytic subset $P$. Then we have
\begin{eqnarray*}
\int_{CC(\One_P)}\ome=\int_{CC(\One_{\bar
P})}\ome-\int_{CC(\One_{\bar P\backslash
P})}\ome=-\int_{CC(\One_{\bar P\backslash P})}\ome.
\end{eqnarray*}
Since $\dim(\bar P\backslash P)<\dim P$ by Proposition
\ref{sg5}(ii) we can use the induction on the dimension of $P$.
Lemma is proved. \qed
\begin{remark}
The differential operator $D$ was introduced by Rumin \cite{rumin}
for an arbitrary contact manifold, and it depends only on the
contact structure. In our case for any smooth manifold $X$ the
space $\PP_+(T^*X)$ has a canonical contact structure, and the
operator $D$ used in the proof of Lemma \ref{im10} corresponds to
it.
\end{remark}

\subsection{Comparison of filtrations.}\label{comparison-filtrations}
Let us define on $\cf(X)$ a filtration by codimension of support:
\begin{eqnarray}\label{ff1}
W_i(\cf(X)):=\{f\in \cf(X)|\, codim (\supp f)\geq i\}.
\end{eqnarray}
We have
$$\cf(X)=W_0(\cf(X))\supset W_1(\cf(X))\supset\dots\supset
W_n(\cf(X))\supset W_{n+1}(\cf(X))=0.$$
\begin{proposition}\label{ff2}
The canonical map
$$\Xi\colon \cf(X)\to \vmi$$
is injective. Moreover for any $i=0,1,\dots,n$, and any $f\in
W_i(\cf(X))\backslash W_{i+1}(\cf(X))$ there exists $\phi\in
W_{n-i}(V^\infty_c(X))$ such that
$$<\Xi(f),\phi>\ne 0.$$
\end{proposition}
{\bf Proof.} Clearly it is enough to prove the second statement.
Let us fix a constructible function $f\in W_i(\cf(X))\backslash
W_{i+1}(\cf(X))$. Thus $\supp f$ is a subanalytic set and
$codim(\supp f)=i$.

One can choose a regular point $x\in \supp f$, a neighborhood $U$,
a real analytic diffeomorphism $g\colon U\tilde\to \RR^n$ such
that $f|_U\circ g^{-1}= c\cdot \One_{\RR^{n-k}}$ where
$\RR^{n-k}\subset \RR^n$ is the coordinate subspace, and $c\ne 0$
is a constant. Thus we may assume that $X=\RR^n,\,
f=\One_{\RR^{n-k}}$. Let us choose $\ome\in
C^\infty_c(T^*\RR^n,\Ome^n\otimes p^*o)$ as follows. Let
$\{(q_1,\dots,q_n)\} $ be coordinates on $\RR^n$. Let
$\{p_1,\dots,p_n\}$ be dual coordinates on $\RR^{n*}$. Let us fix
a $C^\infty$-smooth non-negative compactly supported function
$\tau\colon \RR^{n*}\to \RR_{\geq 0}$ such that $\tau(0)>0$. Let
us take
$$\ome:=\tau\cdot dx_1\wedge\dots \wedge
dx_{n-k}\wedge dy_{n-k+1}\wedge\dots\wedge dy_n.$$ Then clearly
$\int_{CC(\RR^{n-k})}\ome\ne 0$ and $\ome\in
C^\infty_c(T^*\RR^n,W_{n-i}(\Ome^n\otimes p^*o))$. \qed

From now on we will identify $\cf(X)$ with the subspace of
$V^{-\infty}(X)$ via the imbedding $\Xi$.

\begin{proposition}\label{ff3}
(i) $\cf(X)$ is dense in $V^{-\infty}(X)$ in the weak topology.

(ii) For any $i=0,1,\dots,n$
$$\cf(X)\cap W_i(\vmi)=W_i(\cf(X)).$$
\end{proposition}
{\bf Proof.} (i) By the Hahn-Banach theorem it is enough to prove
that for any $\phi\in \vic\backslash \{0\}$ there exists $f\in
\cf(X)$ such that $<f,\phi>\ne 0$. Let us fix $\phi\in
\vic\backslash \{0\}$. One may find an open subset $U\subset X$
and a real analytic diffeomorphism $g\colon U\tilde \to \RR^n$
such that $\phi|_U\not\equiv 0$. The smooth valuation
$g_*\phi|_U\in V^\infty(\RR^n)$ does not vanish identically. By
Proposition 2.4.10 from \cite{part2} there exists a convex compact
set $K\in \ck(\RR^n)\cap \cp(\RR^n)$ such that $(g_*\phi)(K)\ne
0$. Since every compact set can be approximated in the Hausdorff
metric by convex compact polytopes, we may assume that $K$ is a
convex compact polytope, and hence a subanalytic set. Hence
$g^{-1}(K)$ is a compact subanalytic subset of $X$. Take
$f:=\One_{g^{-1}(K)}$. Then $$<f,\phi>\ne 0.$$ Part (i) is proved.

(ii) First let us show the inclusion
\begin{eqnarray}\label{ff4}
\cf(X)\cap W_i(\vmi)\subset W_i(\cf(X)).
\end{eqnarray}
Let $f\in \cf(X)$ be such that $f\not\in W_i(\cf(X))$. Let us
choose $l<i$ such that $f\in W_l(\cf(X))\backslash
W_{l+1}(\cf(X))$. By Proposition \ref{ff2} there exists $\phi\in
W_{n-l}(\vic)$ such that $<f,\phi>\ne 0.$ Hence $f\not\in
W_{l+1}(\vmi)$. Since $l+1\leq i$ we have $f\not\in W_{i}(\vmi)$.
This proves the inclusion (\ref{ff4}).

Let us prove the opposite inclusion
\begin{eqnarray}\label{ff5}
W_i(\cf(X))\subset \cf(X)\cap W_i(\vmi).
\end{eqnarray}
By Proposition \ref{f2} it is enough to show that for any $f\in
W_i(\cf(X)),\, \phi\in W_{n-i+1}(\vic)$
$$<f,\phi>=0.$$ By Lemma \ref{ccc4} and \ref{im10} there exists a {\itshape compactly
supported} form $\ome\in C^\infty_c(T^*X, W_{n-i+1}(\Ome^n\otimes
p^*o))$ such that for any $h\in \cf_c(X)$
\begin{eqnarray}\label{ff5.1}
<h,\phi>=\int_{CC(h)}\ome.
\end{eqnarray}
Since the form $\ome$ is compactly supported the equality
(\ref{ff5.1}) holds for any $h\in \cf(X)$.

Let us assume now that $f=\One_Q$ where $Q$ is a subanalytic
subset with $codim Q\geq i$. We may assume that $Q$ relatively
compact. We have to show that $\int_{CC(Q)}\ome=0$. It is enough
to show that the restriction of $\ome$ to $\supp (CC(Q))$
vanishes. By Lemma \ref{lemma-charact} one can find a finite
covering $\bar Q=\cup_\alp Q_\alp$ such that
$CC(Q)\subset\cup_\alp T^*_{Q_\alp}X$. But since $codim Q_\alp\geq
i$ it is obvious that the restriction of $\ome$ to $T^*_{Q_\alp}X$
vanishes. The proposition is proved. \qed

\subsection{The
integration functional vs. the integration with respect to the
Euler characteristic.}\label{comparison-integration} On the space
$\cf_c(X)$ we have the linear functional $\cf_c(X)\to \CC$ of
integration with respect to the Euler characteristic which is
uniquely characterized by the property $\One_Q\mapsto \chi(Q)$ for
any compact subanalytic subset $Q$ (see \cite{kashiwara-schapira},
\S9.7). For a function $f\in \cf_c(X)$ we will denote the integral
of $f$ with respect to the Euler characteristic by $\int fd\chi$.

Thus we have the canonical imbedding
$$\cf_c(X)\hookrightarrow V_c^{-\infty}(X).$$
\begin{proposition}\label{gg1}
The restriction of the integration functional $\int\colon
V_c^{-\infty}\to \CC$ to $\cf_c(X)$ is equal to the integration
with respect to the Euler characteristic.
\end{proposition}
{\bf Proof.} Since the integration functional $\int\colon \vmic\to
\CC$ to $\cf_c(X)$ is continuous in the weak topology, Proposition
\ref{f10} implies that there exists unique $\xi\in V^\infty(X)$
such that for any $\psi\in \vmic$
$$\int\psi=<\psi,\xi>.$$
It is clear that if $\psi\in \vic$ then
$$\int\psi=<\psi,\chi>.$$
Since $\vic$ is dense in $\vmic$ by Proposition \ref{f5}, it
follows that $\xi=\chi$.

Let us fix a Riemannian metric on $X$. By Theorems 1.5, 1.8 of
\cite{fu-94} there exists a form $\ome\in \tilde
C^\infty(T^*X,\Ome^n\otimes p^*o)$ (which is a little modification
of the Chern-Gauss-Bonnet form \cite{chern}) such that for any
compact subset $P\subset X$ which is either subanalytic or belongs
to $\cp(X)$ one has
$$\chi(P)=\int_{CC(P)}\ome.$$

Then by the construction of the imbedding $\cf(X)\hookrightarrow
\vmi$ and by Proposition \ref{sg6}(ii) we have for any $f\in
\cf_c(X)$
$$\int f=<f,\chi>=\int_{CC(f)}\ome.$$
The proposition is proved. \qed
\subsection{The Euler-Verdier involution vs. the Verdier
duality.}\label{ev-verdier} The space of constructible functions
$\cf(X)$ has a canonical operator
$$\DD\colon\cf(X)\to \cf(X)$$
called the Verdier duality (see \cite{kashiwara-schapira}, \S
9.7). It satisfies $\DD^2=Id$, and for any function $f\in \cf(X)$
\begin{eqnarray}\label{vd1}
CC(\DD f)=a^*CC(f)
\end{eqnarray}
where $a\colon T^*X\to T^*X$ is the antipodal involution
(Proposition 9.4.4 of \cite{kashiwara-schapira}). The main result
of this subsection is the following proposition.
\begin{proposition}\label{vd2}
The restriction of the Euler-Verdier involution $\sigma\colon
\vmi\to \vmi$ to $\cf(X)$ is equal to $(-1)^n\DD$.
\end{proposition}
{\bf Proof.} Let $f\in\cf(X)$. We have to show that for any
$\phi\in \vic$ one has
$$<\sigma f,\phi>=(-1)^n<\DD f,\phi>.$$
By Lemma \ref{ccc4} there exists $\ome\in
C^\infty_c(T^*X,\Ome^n\otimes p^*o)$ such that for any $h\in
\cf(X)$ one has
$$<h,\phi>=\int_{CC(h)}\ome.$$
Then by the definition of $\sigma$ on $V^{-\infty}(X)$ we get
\begin{eqnarray*}
<\sigma f,\phi>=(-1)^n<f,\sigma\phi>=\int_{CC(f)}a^*\ome=\\
(-1)^n\int_{CC(\DD f)}\ome =<(-1)^n\DD f,\phi>.
\end{eqnarray*}
The proposition is proved. \qed

\end{document}